\documentclass[twoside]{amsart}

\usepackage{latexsym, amsmath, amssymb, amsfonts, amscd, amsthm}

\newcommand{\cA}{{\mathcal A}}
\newcommand{\cB}{{\mathcal B}}
\newcommand{\cC}{{\mathcal C}}
\newcommand{\cD}{{\mathcal D}}

\newcommand{\De}{{\mathbb D}}

\newcommand{\diam}{{\rm diam}}
\newcommand{\dist}{{\rm dist}}
\newcommand{\GH}{{\rm GH}}
\newcommand{\CCQM}{{\rm C^*CQM}}
\newcommand{\CM}{{\rm CM}}
\newcommand{\CQM}{{\rm CQM}}
\newcommand{\cq}{{\rm cq}}
\newcommand{\rH}{{\rm H}}
\newcommand{\oq}{{\rm oq}}
\newcommand{\q}{{\rm q}}
\newcommand{\sa}{{\rm sa}}
\newcommand{\im}{{\rm im}}
\newcommand{\Cov}{{\rm Cov}}
\newcommand{\rnu}{{\rm nu}}
\newcommand{\SUB}{{\rm SUB}}
\newcommand{\hol}{{\rm hol}}
\newcommand{\Af}{{\rm Af}}
\newcommand{\mul}{{\rm mul}}
\newcommand{\ru}{{\rm u}}

\newcommand{\Ce}{{\mathbb C}}
\renewcommand{\Re}{{\mathbb R}}

\newcommand{\Ne}{{\mathbb N}}

\newcommand{\pa}{\|}

\newcommand{\ie}{{\it i.e. }}

\theoremstyle{definition}
\newtheorem{theorem}{Theorem}[section]
\newtheorem{corollary}[theorem]{Corollary}
\newtheorem{lemma}[theorem]{Lemma}
\newtheorem{proposition}[theorem]{Proposition}
\newtheorem{remark}[theorem]{Remark}
\newtheorem{notation}[theorem]{Notation}
\newtheorem{example}[theorem]{Example}

\newtheorem{definition}[theorem]{Definition}

\newtheorem*{acknowledgments}{Acknowledgments}

\begin{document}

\title{$C^*$-algebraic Quantum Gromov-Hausdorff Distance}
\author{Hanfeng Li}

\address{Department of Mathematics \\
University of Toronto \\
Toronto ON M5S 3G3, CANADA} \email{hli@fields.toronto.edu}
\date{November 13, 2004}


\begin{abstract}
We introduce a new
quantum Gromov-Hausdorff distance between $C^*$-algebraic
compact quantum metric spaces.
Because it is able to distinguish algebraic structures, this new distance
fixes a weakness of Rieffel's quantum distance.
We show that this new quantum distance has properties analogous to
the basic properties of the classical Gromov-Hausdorff distance, and
we give criteria for when a parameterized
family of $C^*$-algebraic compact quantum metric spaces is continuous
with respect to this new distance.
\end{abstract}

\maketitle

\section{Introduction}
\label{intro:sec}

In  \cite{Connes89} Connes initiated
the study of metric spaces in noncommutative setting
in the framework of his spectral triple \cite{Connes94}.
The main ingredient of a spectral triple is a Dirac operator $D$.  On the
one hand, it captures the differential structure by setting $df=[D,
f]$. On the other hand, it enables us to recover the Lipschitz
seminorm $L$, which is usually defined as
\begin{eqnarray} \label{dist to Lip:eq}
L(f):=\sup \{\frac{|f(x)-f(y)|}{\rho(x, y)}:x \neq y\},
\end{eqnarray}
where $\rho$ is the geodesic metric on the Riemannian manifold,
instead by means of $L(f)=\pa [D, f]\pa$, and then one recovers the metric $\rho$  by
\begin{eqnarray} \label{Lip to dist:eq}
\rho(x, y)=\sup_{L(f)\le 1} |f(x)-f(y)|.
\end{eqnarray}
In Section 2 of \cite{Connes89} Connes went further by considering
the (possibly $+\infty$-valued) metric on the state space of the
algebra defined by (\ref{Lip to dist:eq}). Motivated by what
happens to ordinary compact metric spaces, in \cite{Rieffel98b,
Rieffel99b, Rieffel00} Rieffel introduced ``compact quantum metric
spaces'' which requires the metric on the state space to induce
the weak-$*$ topology. Many interesting examples of compact quantum
metric spaces have been constructed
 \cite{Rieffel98b, Rieffel02, Rieffel03O, Li9}.

Motivated by questions in string theory,
Rieffel also introduced a notion of quantum Gromov-Hausdorff distance
for compact quantum spaces \cite{Rieffel00, Rieffel01}. This is
defined as a
modified ordinary Gromov-Hausdorff distance for the state-spaces.
 The success of this quantum distance is that
there is also a quantum version of the Gromov completeness and
compactness theorems,
which assert that the set of isometry classes of compact metric spaces
endowed with the
Gromov-Hausdorff  distance is
complete, and provide a criterion for when a subset of this complete
space is precompact. And this quantum distance extends the ordinary
Gromov-Hausdorff distance in the sense that sending each compact
metric space to the associated compact quantum metric space is a
homeomorphism (though not an isometry) from the space of isometry
classes of compact metric spaces to a closed subspace of the space of
isometry classes of compact quantum metric spaces.
Since his construction does not involve the
multiplicative structure of the algebras but only needs the
state-spaces, Rieffel set up everything on more general spaces, namely
\emph{order-unit spaces},  based on Kadison's representation theory
\cite{Kadison51, Alfsen71}.
As a consequence, his quantum Gromov-Hausdorff distance can not distinguish
the multiplicative structures of the algebras, \ie
non-isomorphic $C^*$-algebras might have distance zero.

In view of the principle of noncommutative geometry, it may be
more natural to define the quantum distance as a modified
 Gromov-Hausdorff distance for the order-unit spaces (or
 $C^*$-algebras) directly. Under this guidance, in \cite{Li10}
we introduced
an order-unit quantum Gromov-Hausdorff distance as a modified ordinary
Gromov-Hausdorff distance for certain balls in the order-unit spaces.
It turns out that this order-unit quantum distance is Lipschitz equivalent to
Rieffel's quantum distance \cite[Theorem 1.1]{Li10}.
However, as an advantage of our approach,
we gave criteria for when a
parameterized family of compact quantum metric spaces is continuous with
respect to the order-unit quantum distance \cite[Theorems 1.2 and 7.1]{Li10}.

Another advantage of our approach is that we can adapt it easily to
deal with additional algebraic structures.
This is what we shall do in this paper.
We introduce
a \emph{$C^*$-algebraic quantum Gromov-Hausdorff distance}, and show that it distinguishes
the multiplicative structures of the algebras. We also show that
it  extends the ordinary
Gromov-Hausdorff distance in the same way as does Rieffel's quantum
distance. There is also a quantum version of the Gromov completeness and
compactness theorems, though in order to make sure that the limits are
$C^*$-algebras we need some compatibility condition between the
seminorms and the algebraic structures, which turns out to be some kind of
generalized Leibniz rule. We also give similar criteria for when a
parameterized family of $C^*$-algebraic compact quantum metric spaces is continuous with
respect to the $C^*$-algebraic quantum distance.
These criteria enable us to conclude
immediately that the continuity of noncommutative
tori \cite[Theorem 9.2]{Rieffel00} and $\theta$-deformations \cite[Theorem 1.4]{Li10}
and the convergence of matrix algebras to integral coadjoint orbits of compact connected semisimple
Lie groups \cite[Theorem 3.2]{Rieffel01} all hold with respect to the $C^*$-algebraic quantum distance.

 David Kerr \cite{Kerr02} has defined a matricial quantum
distance which is also
able to distinguish the multiplicative structures of algebras.
His approach is based on operator systems, very different from ours.
In \cite{KL} we apply the results and methods of this paper to show that
this distance has similar properties as our $C^*$-algebraic quantum distance does.

This paper is organized as follows.
In Section~\ref{Prelim:sec} we review briefly the Gromov-Hausdorff distance for
ordinary compact metric spaces
 and Rieffel's theory of compact quantum metric spaces.
The $C^*$-algebraic quantum Gromov-Hausdorff distance
is introduced in Section~\ref{C*QGH:sec}.
We show there that this distance distinguishes multiplicative structures
of the algebras.
Section~\ref{C*QGCC:sec} is devoted to a quantum version
of Gromov completeness and compactness theorems.
The criteria for when a
parameterized family of $C^*$-algebraic compact quantum metric spaces is continuous with
respect to the $C^*$-algebraic quantum distance are given in Section~\ref{COntFieldC*QCM:sec}.

\begin{acknowledgments}
  This is part of my Ph.D. dissertation submitted to UC Berkeley in 2002.
  I am indebted to my advisor, Professor Marc Rieffel, for many helpful discussions, suggestions,
  and for his support throughout my time at Berkeley.
I am grateful to Nik Weaver for providing me Example~\ref{modified 4 point:eg}.
I thank David Kerr for sending me his preprints and helpful discussions.
I also thank Thomas Hadfield and Fr\'ed\'eric Latr\'emoli\`ere
  for valuable conversations.
\end{acknowledgments}

\section{Preliminaries}
\label{Prelim:sec}

In this section we review briefly the Gromov-Hausdorff distance for
compact metric spaces \cite{Gromov99,
Sakai96, BBI01} and Rieffel's theory of compact quantum metric spaces
\cite{Rieffel99b, Rieffel00, Rieffel03}.

Let $(X, \rho)$ be a metric space, \ie $\rho$ is a metric on the space
$X$. For any subset $Y\subseteq X$ and $r>0$, let
\begin{eqnarray*}
\cB_r(Y)=\{x\in X:\rho(x,y)<r \mbox{ for some } y\in Y\}
\end{eqnarray*}
be the set of points with distance less than $r$ from $Y$.
When $Y=\{x\}$, we also write it as
$\cB_r(x)$ and call it the {\it open ball} of radius $r$
centered at $x$.

For nonempty subsets $Y, Z\subseteq X$, we can measure the
distance between $Y$ and $Z$ inside of $X$ by the \emph{Hausdorff distance}
$\dist^{\rho}_{\rH}(Y, Z)$ defined by
\begin{eqnarray*}
 \dist^{\rho}_{\rH}(Y, Z):=\inf \{r>0:Y\subseteq \cB_r(Z),
 Z\subseteq \cB_r(Y)\}.
\end{eqnarray*}
We will also use the notation $\dist^X_{\rH}(Y, Z)$ when there is no
confusion about the metric on $X$. When $X$ is compact, the set $\SUB(X)$ of nonempty closed subsets
of $X$ is compact equipped with the metric $\dist^X_{\rH}$ \cite[Proposition 7.3.7]{BBI01}.

For any compact metric spaces $X$ and $Y$, Gromov \cite{Gromov81}
introduced the \emph{Gromov-Hausdorff
 distance}, $\dist_{\GH}(X, Y)$, which is defined as
\begin{eqnarray*}
\dist_{\GH}(X, Y)&:=&\inf \{\dist^{Z}_{\rH}(h_X(X), h_Y(Y))|\,
h_X:X\rightarrow Z, h_Y:Y\rightarrow Z \mbox{ are }\\ & & \mbox{
isometric embeddings into some metric space } Z\}.
\end{eqnarray*}
It is possible to reduce the space $Z$ in above to be the disjoint union $X\coprod Y$. A
distance $\rho$ on $X\coprod Y$ is said to be {\it admissible} if the
inclusions $X, Y\hookrightarrow X\coprod Y$ are isometric
embeddings. Then it is not difficult to check that
\begin{eqnarray*}
\dist_{\GH}(X, Y)=\inf \{ \dist^{\rho}_{\rH}(X, Y):\rho
\mbox{ is an admissible distance on } X\coprod Y\}.
\end{eqnarray*}

For a compact metric space $(X, \rho)$, we shall denote by
$\diam(X):=\max \{\rho(x, y)|\\x, y\in X\}$
the \emph{diameter} of $X$.
Also let $r_X=\frac{\diam(X)}{2}$ be the \emph{radius} of $X$.
For any $\varepsilon>0$, the \emph{covering number}
$\Cov_{\rho}(X, \varepsilon)$ is defined as the smallest number of open balls of
radius $\varepsilon$ whose union covers $X$.

Denote by $\CM$ the set of isometry classes of compact metric spaces.
One important property of Gromov-Hausdorff distance is the completeness and
compactness theorem by Gromov \cite{Gromov81}:

\begin{theorem}[Gromov's Completeness and Compactness Theorems] \label{GH:thm}
The space $(\CM, \dist_{\GH})$ is a complete metric space. A subset
$\mathcal{S}\subseteq \CM$ is totally bounded
(\ie has compact closure) if and only if

(1) there is a constant $D$ such that $\diam(X, \rho)\le D$ for all
$(X, \rho)\in \mathcal{S}$;

(2) for any $\varepsilon>0$, there exists a constant $K_{\varepsilon}>0$ such that
$\Cov_{\rho}(X, \varepsilon)\le K_{\varepsilon}$ for all $(X, \rho)\in \mathcal{S}$.

\end{theorem}

Next we recall Rieffel's definition of ($C^*$-algebraic) compact quantum
metric spaces.

Though Rieffel has set up his theory in the general framework of
order-unit spaces, we shall  mainly need it only for $C^*$-algebras.
See the discussion preceding Definition 2.1 in \cite{Rieffel00}
for the reason of requiring the reality condition (\ref{real:eq}) below.

\begin{definition}\cite[Definition 2.1]{Rieffel00} \label{C*QCM:def}
By a \emph{$C^*$-algebraic compact quantum metric
space} we mean a pair $(\cA, L)$ consisting of a unital
$C^*$-algebra $\cA$ with identify $e_{\cA}$ and  a (possibly
$+\infty$-valued) seminorm $L$ on $\cA$ satisfying
 the \emph{reality condition}
\begin{eqnarray} \label{real:eq}
L(a)&=& L(a^*),
\end{eqnarray}
such that
$L$ vanishes on $\Ce e_{\cA}$ and the metric $\rho_L$ on the state space
$S(\cA)$ defined by (\ref{Lip to dist:eq}) induces the weak-$*$ topology.
The \emph{diameter}, \emph{radius}  and \emph{covering number}
of $(\cA, L)$ are defined to be the same as those of $(S(\cA), \rho_L)$.
We say that $L$ is a \emph{Lip-norm}.
\end{definition}

Let $\cA$ be a unital $C^*$-algebra and let $L$ be a
(possibly
$+\infty$-valued) seminorm $L$ on $\cA$ satisfying
 (\ref{real:eq}).
Then one can check easily that the pair $(\cA, L)$ is a
$C^*$-algebraic compact quantum metric space if and only if the
space $A=\{a\in \cA_{\sa}: L(a)<\infty\}$ is dense in
$\cA_{\sa}$ containing $e_{\cA}$ and the pair $(A,
L|_A)$ is a compact quantum metric space in the sense of
\cite[Definition 2.1]{Rieffel00}.
We shall call the latter the associated compact quantum metric space of $(\cA, L)$.
For most parts of this paper the reader may take the compact quantum metric
spaces as the ones associated to $C^*$-algebraic compact quantum metric spaces.

For any compact metric space $(X, \rho)$, define the Lipschitz
seminorm (which may take value $+\infty$), $L_\rho$, on $C_{\Ce}(X)$,
the set of all $\Ce$-valued continuous functions over $X$, by
(\ref{dist to Lip:eq}).
Then the pair $(C_{\Ce}(X), L_{\rho})$ is a $C^*$-algebraic
compact quantum metric space, and the associated compact quantum metric space
is the pair $(A_X, L_{\rho}|_{A_X})$, where $A_X$ is the space of $\Re$-valued
Lipschitz functions on $X$.
They are called the associated
  ($C^*$-algebraic) compact quantum
metric space of $(X, \rho)$.
Denote by $\CQM$
the set of isometry classes of compact quantum metric spaces \cite[Definition 6.3]{Rieffel00}.
In \cite[Theorem 13.16]{Rieffel00}
Rieffel showed that the map $(X, \rho)\mapsto (A_X, L|_{A_X})$ from
$(\CM, \dist_{\GH})$ to $(\CQM, \dist_{\q})$ is continuous,
injective, and sends $ \CM$ to a closed subspace of $\CQM$, where
$\dist_{\q}$ is Rieffel's quantum distance. The argument there actually shows that
this map is a homeomorphism from $ \CM$ to its image.


\section{Definition of $C^*$-algebraic quantum distance}
\label{C*QGH:sec}

In this section we define a quantum Gromov-Hausdorff distance $\dist_{\cq}$
which distinguishes the algebraic structures. The main results in this section are
Theorems~\ref{basic C*QGH:thm} and \ref{C*GH GH:thm}.

\begin{notation} \label{ball:notation}
Let $(\cA, L)$ be a $C^*$-algebraic compact quantum metric space.
For any $r\ge 0$, set
\begin{eqnarray*}
\cD_r(A):=\{a\in A: L(a)\le 1, \pa a\pa \le r\}, \quad \mbox{ and } \quad \cD(A):=\cD_{r_{\cA}}(A),
\end{eqnarray*}
where $(A, L|_A)$ is the compact quantum metric space associated to $(\cA, L)$
(see the discussion after Definition~\ref{C*QCM:def}).
$\cD_r(\cA)$ and $\cD(\cA)$ are defined similarly.
\end{notation}

We require $\dist_{\cq}$ to extend the ordinary Gromov-Hausdorff
distance $\dist_{\GH}$ in the sense that they define the same
topology on $\CM$ (see the discussion at the end of
Section~\ref{Prelim:sec}). To make $\dist_{\cq}$ manageable, we also
want an analogue of the criteria for continuity of parameterized
families of compact quantum metric spaces with respect to the
order-unit quantum distance $\dist_{\oq}$ \cite[Theorems 1.2 and 7.1]{Li10}. In
\cite[Definition 4.2]{Li10} the order-unit distance $\dist_{\oq}$ is
defined as a modified Gromov-Hausdorff distance between the balls
$\cD(A)$, based on the observation that these
balls catch both the norm and the Lip-norm of the compact quantum
metric spaces \cite[Lemma 4.1]{Li10}. Precisely,
for any $R\ge 0$ and
any compact quantum metric spaces $(A, L_A)$ and $(B, L_B)$ we have
\begin{eqnarray*}
\dist_{\oq}(A, B)=\inf \{\max (\dist^V_{\rH}(h_A(\cD(A)), h_B(\cD(B))),
\pa h_A(r_Ae_A)-h_B(r_Be_B)\pa)\},\\
\dist^R_{\oq}(A, B) = \inf \{\max
(\dist^V_{\rH}(h_A(\cD_{ R}(A)), h_B(\cD_{ R}(B))),
\pa h_A(Re_A)-h_B(Re_B)\pa) \},
\end{eqnarray*}
where the infima are taken over all triples $(V, h_A, h_B)$
consisting of a real normed space $V$ and linear isometric
embeddings $h_A:A\rightarrow V$ and $h_B:B\rightarrow V$. The
distance $\dist^R_{oq}$ was introduced to facilitate certain
arguments involving continuity.

If we look at above definition of $\dist_{\oq}$,
one natural choice is to define
$\dist_{\cq}(\cA, \cB)$ as
$\inf\{\dist^{\cC}_{\rH}(h_{\cA}(\cD(\cA)),
h_{\cB}(\cD(\cB)))\}$, where the infimum is taken over all
faithful $*$-homomorphisms $h_{\cA}$ and $h_{\cB}$ of
$\cA$ and $\cB$ into some
$C^*$-algebra $\cC$. There are two difficulties. The first one
is that $\Ce$-valued functions do not behave as well as $\Re$-valued
functions do under Lipschitz seminorm. So if we use the balls $\cD(\cA)$,
the resulting quantum distance will
induce a stronger topology on $\CM$ than $\dist_{\GH}$ does
(see Remark~\ref{Weaver:remark}). The second difficulty is that we
do not have the $C^*$-algebraic analogue of \cite[Lemma 7.2]{Li10}
(see \cite[Remark 7.5]{Li10}), which is crucial in the proofs of
\cite[Theorems 1.2 and 7.1]{Li10}. Hence if we require
$h_{\cA}$ to be a $*$-homomorphism, then we do not know how
to prove the  analogues of \cite[Theorems 1.2 and 7.1]{Li10} (however,
see Remark~\ref{nuclear:remark}). Thus we are forced to use the
ball $\cD(A)$ and to require
$h_{\cA}$ to be a linear isometric embedding into a
complex normed space $V$. But this will not distinguish
the algebraic structures. To get around this difficulty, we will use the
graph of the multiplication in $\cD(A)$
instead of $\cD(A)$ itself.

\begin{notation} \label{graph:notation}
For a subset $X$ of a $C^*$-algebra $\cA$ we shall denote by
\begin{eqnarray*}
X^m:=\{(x,y,xy)\in \cA\oplus \cA\oplus \cA:x,
y\in X\}
\end{eqnarray*}
the graph of the multiplication in $X$. For a normed space $V$ we shall
denote by $V^{(3)}$ the direct sum $V\oplus V\oplus V$ equipped with
the supremum-norm $\pa (v_1, v_2, v_3)\pa=\max(\pa v_1\pa, \pa v_2\pa, \pa
v_3\pa) $. Similarly, we define $V^{(2)}$.
For a linear map $h:\cA\to V$ we denote by
$h^{(3)}$ the natural extension $\cA^{(3)}\to V^{(3)}$.
\end{notation}

In the above definition of $\dist_{oq}$
we had to consider
the term $\pa h_A(r_Ae_A)-h_B(r_Be_A)\pa$ to keep track
of the units.
Since the algebraic structures force
$\pa h_{\cA}(e_{\cA})-h_{\cB}(e_{\cB})\pa $
to be bounded by the Hausdorff distance between
$h^{(3)}_{\cA}((\cD(A))^m)$ and
$h^{(3)}_{\cB}((\cD(B))^m)$ (see
Lemma~\ref{cq o-u:lemma}),
we do not need this term
any longer as we did in $\dist_{\oq}$.

\begin{definition} \label{C*QGH:def}
Let $(\cA, L_{\cA})$ and $(\cB,L_{\cB})$ be
$C^*$-algebraic compact quantum metric
spaces. We define the \emph{$C^*$-algebraic quantum Gromov-Hausdorff distance} between them,
denoted by $\dist_{\cq}(\cA, \cB)$, by
\begin{eqnarray*}
\dist_{\cq}(\cA, \cB):=
\inf \{\dist^{V^{(3)}}_{\rH}(h^{(3)}_{\cA}((\cD(A)^m)),
h^{(3)}_{\cB}((\cD(B)^m))\},
\end{eqnarray*}
and, for $R\ge 0$,  the \emph{$R$-$C^*$-algebraic quantum
Gromov-Hausdorff distance} between them, denoted by
$\dist^R_{\cq}(\cA, \cB)$, by
\begin{eqnarray*}
\dist^R_{\cq}(\cA, \cB):=
\inf \{\dist^{V^{(3)}}_{\rH}(h^{(3)}_{\cA}((\cD_R(A)^m)),
h^{(3)}_{\cB}((\cD_R(B)^m))\},
\end{eqnarray*}
where the infima are taken over all triples $(V, h_{\cA}, h_{\cB})$ consisting of a complex normed space
$V$ and linear isometric embeddings $h_{\cA}:\cA\rightarrow V$ and $h_{\cB}:\cB\rightarrow V$.
\end{definition}

\begin{remark} \label{Kadison:remark}

(1) To simply the notation, usually we shall identify $\cA$ and $\cB$ with
    their images $h_{\cA}(\cA)$ and $h_{\cB}(\cB)$ respectively,
    and just say that
    $V$ is a normed space containing both $\cA$ and $\cB$;

(2) Initiated by Kadison and Kastler \cite{KS72} there has been a lot of work
on perturbation of $C^*$-algebras \cite{CC83, PR81, Johnson94}. In this theory people consider
$\dist^{\cC}_{\rH}(\De(\cA),
  \De(\cB))$, where $\cA$ and $\cB$ are
  $C^*$-subalgebras of a $C^*$-algebra $\cC$, and
  $\De(\cA)$ is the unit ball of $\cA$.
Using Blackadar's result \cite[Theorem 3.1]{Blackadar80} on amalgamation of $C^*$-algebras one sees
easily that $\inf(\dist^{\cC}_{\rH}(\De(\cA),
  \De(\cB)))$ defines a distance $\dist_{\ru}$, where $\cC$
  runs over all $C^*$-algebras containing both $\cA$ and
  $\cB$ as $C^*$-subalgebras. We do not know when $\dist_{\ru}(\cA,
  \cB)=0$.
  This looks very similar to the distance $\dist_{\rnu}$
  we shall discuss in Remark~\ref{nuclear:remark}. However, they are
  quite different. If $\cA$ is a separable unital continuous trace $C^*$-algebra,
  then when $\dist_{\ru}(\cA, \cB)$ is small,
  $\cA$ and $\cB$
  must be isomorphic as $C^*$-algebras \cite[Theorem 4.22]{PR81}.
\end{remark}

For any normed
spaces $V$ and $W$ we call
a norm $\pa \cdot \pa_{V\oplus W}$ on $V\oplus W$ \emph{admissible} if
it extends the norms on $V$ and $W$.
By an argument similar to that in the proof of \cite[Proposition 4.4]{Li10}
we have:

\begin{proposition} \label{V=A+B 2:prop}
Let $(\cA, L_{\cA})$ and $(\cB,L_{\cB})$ be
$C^*$-algebraic compact quantum metric
spaces. Then
\begin{eqnarray*}
\dist_{\cq}(\cA, \cB)=
\inf \{\dist^{(\cA\oplus \cB)^{(3)}}_{\rH}((\cD(A))^m,
(\cD(B))^m)\},
\end{eqnarray*}
and, for any $R\ge 0$,
\begin{eqnarray*}
\dist^R_{\cq}(\cA, \cB)=
\inf \{\dist^{(\cA\oplus \cB)^{(3)}}_{\rH}((\cD_R(A))^m,
(\cD_R(B))^m)\},
\end{eqnarray*}
where the infima are taken over all admissible norms on $\cA\oplus
\cB$.
\end{proposition}

Using the amalgamation of normed spaces \cite[Lemma 4.5]{Li10}
we have the triangle inequalities for
$\dist_{\cq}$ and $\dist^R_{\cq}$:

\begin{lemma} \label{triangle dist_cq:lemma}
For any $C^*$-algebraic compact quantum metric spaces $(\cA, L_{\cA}),
(\cB,L_{\cB})$
and $(\cC,L_{\cC})$ we
have
\begin{eqnarray*}
\dist_{\cq}(\cA, \cC)\le \dist_{\cq}(\cA, \cB)+\dist_{\cq}(\cB, \cC).
\end{eqnarray*}
For $R\ge 0$ we also have
\begin{eqnarray*}
\dist^R_{\cq}(\cA, \cC)\le \dist^R_{\cq}(\cA,
\cB)
+\dist^R_{\cq}(\cB,\cC).
\end{eqnarray*}
\end{lemma}

One can take closures for compact quantum metric spaces
\cite[Definition 4.5]{Rieffel99b} (see also \cite[Section 2]{Li10}).
Similarly, we can take closure for
$C^*$-algebraic compact quantum metric spaces:

\begin{lemma} \label{clo of C*QCM:lemma}
Let $(\cA, L)$ be a $C^*$-algebraic compact quantum metric
space. Define a seminorm $L^c$ on $\cA$ by
\begin{eqnarray*}
L^c(a):=\inf \{\liminf_{n\to \infty}L(a_n):a_n\in \cA, \lim_{n\to
  \infty}a_n=a\}.
\end{eqnarray*}
Then $(\cA, L^c)$ is also a $C^*$-algebraic compact quantum metric
space, and the associated compact quantum metric space is exactly the
closure of $(A, L|_A)$.
\end{lemma}
\begin{proof}Since $L(a^*)=L(a)$ for all $a\in \cA$, clearly
  $L^c(a^*)=L^c(a)$ for all $a\in \cA$ too.

For any $a\in \cA_{\sa}$, if $a_n\to a$ when $n\to \infty$,
then also $(a_n)_{\sa}\to a$. Notice that
$L((a_n)_{\sa})\le \frac{L(a_n)+L(a^*_n)}{2}=L(a_n)$.
So $\liminf_{n\to \infty}L((a_n)_{\sa})\le
\liminf_{n\to \infty}L(a_n)$. Hence
\begin{eqnarray*}
L^c(a)=\inf \{\liminf_{n\to \infty}L(a_n):a_n\in \cA_{\sa},
  \lim_{n\to \infty}a_n=a\}
\end{eqnarray*}
for $a\in \cA_{\sa}$.
It follows immediately that $L^c$ restricted to $\cA_{\sa}$ is exactly
the closure of $L|_A$.
So the set
$\{a\in \cA_{\sa}:L^c(a)<\infty\}$ equipped with the
restriction of $L^c$
is a compact quantum metric space and is actually the closure of $(A,
L|_A)$.
\end{proof}

\begin{definition} \label{clo of C*QCM:def}
Let $(\cA, L)$ be a $C^*$-algebraic compact quantum metric
space. We call $(\cA, L^c)$ the \emph{closure} of $(\cA,
L)$. We say that $(\cA, L)$ and $L$ are \emph{closed} if $L^c=L$.
\end{definition}

Clearly we have $\dist_{\cq}((\cA, L), (\cA, L^c))=
\dist^R_{\cq}((\cA, L), (\cA, L^c))=0$.

Next we establish some basic inequalities about $\dist_{\cq}$ and
$\dist^R_{\cq}$:

\begin{proposition} \label{dist_cq<dist^R_cq:prop}
For any $C^*$-algebraic compact quantum metric spaces $(\cA,
L_{\cA})$ and $(\cB,L_{\cB})$  we
have
\begin{eqnarray}\label{dist_cq 1:eq}
|r_{\cA}-r_{\cB}|\le \dist_{\GH}(\cD(A),\cD(B))
&\le &\dist_{\cq}(\cA, \cB)
\le
\max(r_{\cA}, \, r_{\cB},\, r^2_{\cA}, \, r^2_{\cB}),\\
 \label{dist_cq 2:eq}
|\dist_{\cq}(\cA, \cB)-\dist^{r_B}_{\cq}(\cA,
    \cB)|
&\le& \max(|r_{\cA}-r_{\cB}|,\, |r^2_{\cA}-r^2_{\cB}|).
\end{eqnarray}
For any $R\ge r\ge 0$ we also have
\begin{eqnarray} \label{dist_cq 3:eq}
\dist^r_{\cq}(\cA, \cB)\le \dist^R_{\cq}(\cA,
\cB)\max(2, \,
\dist^R_{\cq}(\cA, \cB)+2r+1).
\end{eqnarray}
\end{proposition}
\begin{proof}
For any compact metric spaces $X$ and $Y$, one has
\begin{eqnarray} \label{GH:eq}
|r_X-r_Y|\le \dist_{\GH}(X, Y).
\end{eqnarray}
The inequality (\ref{dist_cq 1:eq}) follows from (\ref{GH:eq})
and the fact that for any normed space
  $V$ containing $\cA$ and $\cB$
we have
\begin{eqnarray*}
\dist^{V}_{\rH}(\cD(A),
\cD(B)) \le  \dist^{V^{(3)}}_{\rH}((\cD(A))^m,
(\cD(B))^m) \le
 \max(r_{\cA}, \, r_{\cB},\,
 r^2_{\cA}, \, r^2_{\cB}).
\end{eqnarray*}
 To show (\ref{dist_cq 2:eq}) it suffices
to show that $\dist^{\cA^{(3)}}_{\rH}((\cD(A))^m,
(\cD_{r_{\cB}}(A))^m)\le
  \max(|r_{\cA}-r_{\cB}|,\, |r^2_{\cA}-r^2_{\cB}|)$. In fact we have:
\begin{lemma} \label{R, r 2:lemma}
For any $C^*$-algebraic compact quantum metric space $(\cA,
L_{\cA})$ and any $R> r\ge 0$ we have
\begin{eqnarray*}
\dist^{\cA^{(3)}}_{\rH}((\cD_R(A))^m,
(\cD_r(A))^m)\le \max(R-r,\, R^2-r^2).
\end{eqnarray*}
\end{lemma}
\begin{proof}
Notice that $(\cD_r(A))^m$ is
a subset of $(\cD_R(A))^m$.
For any $a_1, a_2\in \cD_R(A)$,
let $a'_j=\frac{r}{R}a_j$. Then $a'_j\in
\cD_r(A)$ and $\pa (a_1, \, a_2, \, a_1a_2)-(a'_1, \, a'_2, \, a'_1a'_2)\pa \le \max(R-r,\, R^2-r^2)$.
This yields the desired inequality.
\end{proof}
We proceed to show (\ref{dist_cq 3:eq}).
We may assume that $r>0$ and both $L_{\cA}$ and $L_{\cB}$ are closed.
Let $V$ be a complex normed
space containing
$\cA$ and $\cB$,  and let
$d=\dist^{V^{(3)}}_{\rH}((\cD_R(A))^m,
(\cD_R(B))^m)$.
For any $a_1, a_2\in \cD_r(A)$ pick $b_1, b_2\in \cD_R(B)$
such that $\pa a_1-b_1\pa, \pa a_2-b_2\pa, \pa
a_1a_2-b_1b_2\pa \le d$.
%
Then $\pa b_j\pa \le d+r$.
Let $\eta=r/(d+r)$, and let
$b'_j=\eta b_j$.
One checks easily that
$b'_j\in  \cD_r(B)$, and
$\pa a_j-b'_j\pa, \, \pa
a_1a_2-b'_1b'_2\pa \le d\cdot \max(2, d+2r+1) $.
Similarly, one can deal with
  pairs in $\cD_r(B)$. Therefore $\dist^{V^{(3)}}_{\rH}((\cD_r(A))^m,
(\cD_r(B))^m)\le d\cdot \max(2, d+2r+1) $.
Then (\ref{dist_cq 3:eq}) follows.
\end{proof}

\begin{lemma} \label{cq o-u:lemma}
Let $(\cA, L_{\cA})$ and
$(\cB,L_{\cB})$ be $C^*$-algebraic compact quantum
metric spaces. Let $V$ be a complex normed space containing $\cA$
and $\cB$. Then for any $R\ge 0$
we have
\begin{eqnarray*}
R^2\pa e_{\cA}-e_{\cB}\pa\le (4R+1)\dist^{V^{(3)}}_{\rH}((\cD_R(A))^m, (\cD_R(B))^m).
\end{eqnarray*}
\end{lemma}
\begin{proof}
Set $d=\dist^{V^{(3)}}_{\rH}((\cD_R(A))^m,
(\cD_R(B))^m)$.
Denote $\cD_R(A)$ and $\cD_R(B)$ by $X$ and $Y$ respectively.
We may assume that both $(\cA, L_{\cA})$ and
$(\cB,L_{\cB})$ are closed. Then $X$ and $Y$
are compact.
Thus for
$(Re_{\cA}, Re_{\cA}, R^2e_{\cA})\in
  X^m$ we can find $(b, b', bb')\in Y^m$ such that $\pa
  Re_{\cA}-b\pa,\, \pa Re_{\cA}-b'\pa,\, \pa
  R^2e_{\cA}-bb'\pa\le d$. Also for
  $(Re_{\cB}, b, Rb)\in
  Y^m$ we can find $(a, a', aa')\in X^m$ such that $\pa
  Re_{\cB}-a\pa,\, \pa b-a'\pa,\, \pa
  Rb-aa'\pa\le d$. Then
\begin{eqnarray*}
\pa a'-Re_{\cA}\pa &\le & \pa a'-b\pa +\pa b-Re_{\cA}\pa
\le 2d, \\
\pa Rb-R^2e_{\cB}\pa &\le &\pa Rb-aa'\pa+\pa aa'-Ra\pa +\pa
Ra-R^2e_{\cB}\pa \\ &\le &
d+\pa a\pa \cdot(\pa a'-b\pa+\pa b-Re_{\cA}\pa)+Rd
\le d+3Rd, \\
\pa R^2e_{\cA}-R^2e_{\cB}\pa &\le&  \pa
R^2e_{\cA}-Rb\pa +\pa Rb-R^2e_{\cB}\pa
\le  Rd+d+3Rd=4Rd+d,
\end{eqnarray*}
as desired.
\end{proof}

\begin{proposition} \label{dist^R_oq<dist^R_cq:prop}
Let $(\cA, L_{\cA})$ and
$(\cB,L_{\cB})$ be $C^*$-algebraic compact quantum
metric spaces. For any $R>0$
we have
\begin{eqnarray}
\dist^R_{\oq}(A, B)\le  \frac{4R+1}{R}\dist^R_{\cq}(\cA,
\cB). \label{dist_cq 4:eq}
\end{eqnarray}
If furthermore $R\ge \max(r_{\cA}, r_{\cB})$, then
\begin{eqnarray}
\frac{2R}{5(4R+1)}|r_{\cA}-r_{\cB}|\le \dist^R_{\cq}(\cA,
\cB). \label{dist_cq 5:eq}
\end{eqnarray}
\end{proposition}
\begin{proof}
(\ref{dist_cq 4:eq}) follows from Lemma~\ref{cq o-u:lemma}
and the definitions of $\dist^R_{\oq}$ and $\dist^R_{\cq}$.
Also (\ref{dist_cq 5:eq}) follows from (\ref{GH:eq}), (\ref{dist_cq 4:eq}),
\cite[Theorem 1.1]{Li10}, and the fact $\dist_{\GH}(S(A), S(B))\le \dist_{\q}(A, B)$.
\end{proof}

\begin{proposition} \label{dist^R_cq<dist^r_cq:prop}
Let $(\cA, L_{\cA})$ and
$(\cB,L_{\cB})$ be $C^*$-algebraic compact quantum
metric spaces. For any $R\ge r\ge \max(r_{\cA},
r_{\cB})>0$
set $\tau_1=(4r+1)(R+r)/(r^2)+1$, and set
$\tau_2=(\tau_1+1)(R+r)+2\tau_1R+1$.
Then we have
\begin{eqnarray} \label{dist_cq 6:eq}
\dist^R_{\cq}(\cA, \cB)\le \dist^r_{\cq}(\cA,
\cB)\max(2\tau_1,
\tau_1\dist^r_{\cq}(\cA, \cB)+\tau_2).
\end{eqnarray}
\end{proposition}
\begin{proof}
Set $X_R=\cD_R(A)$, and set $Y_R=\cD_R(B)$.
Similarly define $X_r$ and $Y_r$.
Let $V$ be a complex normed space containing
$\cA$ and $\cB$. Denote
$\dist^{V^{(3)}}_{\rH}((X_r)^m,(Y_r)^m)$ by $d$.
We may assume that$(\cA, L_{\cA})$ and
$(\cB,L_{\cB})$ are both closed.
By Lemma~\ref{cq o-u:lemma} we have $\pa
  e_{\cA}-e_{\cB}\pa \le ((4r+1)d)/(r^2)$.
Let $a_1,\, a_2 \in X_R$. By \cite[Lemma 4.1]{Li10} we
  can write $a_j$ as $\lambda_j e_{\cA}+a'_j$ with $\lambda_j\in
  \Re$ and $a'_j\in X_r$. Then
$|\lambda_j|\le  R+r$.
Pick $b'_1,\, b'_2 \in Y_r$ such that
$\pa a'_1-a'_2\pa, \, \pa a'_2-b'_2\pa,\,  \pa a'_1a'_2-b'_1b'_2\pa\le
  d$. Set $b_j=\lambda_je_{\cB}+b'_j$.
A routine calculation yields
\begin{eqnarray*}
\pa a_j-b_j\pa \le \tau_1d, \quad \pa b_j\pa \le \tau_1d+R, \quad \pa a_1a_2-b_1b_2\pa\le ((\tau_1+1)(R+r)+1)d.
\end{eqnarray*}
Set $\eta=R/(\tau_1d+R)$, and set
$b''_j=\eta b_j$. Then $b''_j\in Y_R$, and
\begin{eqnarray*}
\pa b_j-b''_j\pa \le \tau_1d,\quad
\pa b_1b_2-b''_1b''_2\pa \le \tau_1^2d^2+2\tau_1dR.
\end{eqnarray*}
Consequently,
\begin{eqnarray*}
\pa a_j-b''_j\pa \le 2\tau_1d, \quad
\pa a_1a_2-b''_1b''_2\pa \le \tau_1^2d^2+\tau_2d.
\end{eqnarray*}
Similarly, we can deal with pairs in $Y_R$. Thus
\begin{eqnarray*}
\dist^{V^{(3)}}_{\rH}((X_R)^m,
(Y_R)^m)\le d\cdot \max(2,\tau_1^2d+\tau_2).
\end{eqnarray*}
Letting $V$ run over all complex normed spaces containing
 $\cA$ and $\cB$, we get (\ref{dist_cq 6:eq}).
\end{proof}

To show that $\dist_{\cq}$ is a metric, we need to see what happens when
$\dist_{\cq}(\cA, \cB)=0$. As promised at the beginning
of this section, isometric $C^*$-algebraic compact quantum
metric spaces are also  isomorphic as $C^*$-algebras:

\begin{definition} \label{iso of C*QCM:def}
Let $(\cA, L_{\cA})$ and $(\cB,
L_{\cB})$ be $C^*$-algebraic compact quantum metric
spaces. By an {\it isometry} from $(\cA,
L_{\cA})$ to $(\cB, L_{\cB})$ we mean a
$*$-isomorphism $\varphi$ from $\cA$ onto
$\cB$ such that $L^c_{\cA}=L^c_{\cB}\circ
\varphi$ on $\cA_{\sa}$. When there exists an isometry from
$(\cA, L_{\cA})$ to $(\cB,
L_{\cB})$, we say that $(\cA, L_{\cA})$ and
$(\cB, L_{\cB})$ are \emph{isometric}.
We denote by $\CCQM$ the set of all isometry classes of
$C^*$-algebraic compact quantum metric spaces, and denote by
$\CCQM^R$ the subset consisting of isometry classes of
$C^*$-algebraic compact quantum metric spaces with radii no bigger
than $R$.
\end{definition}

Notice that we require $L^c_{\cA}=L^c_{\cB}\circ
\varphi$ only on $\cA_{\sa}$, not on the whole of
$\cA$.

\begin{theorem} \label{basic C*QGH:thm}
$\dist_{\cq}$ and $\dist^R_{\cq}$ are metrics on $\CCQM$ and
$\CCQM^R$ respectively.
They define the same topology on
$\CCQM^R$. The forgetful map
$(\cA, L)\mapsto (A, L|_A)$ from
$(\CCQM, \dist_{\cq})$ to
$(\CQM, \dist_{\oq})$
is
continuous.
\end{theorem}
\begin{proof} It suffices to show that $(\cA, L_{\cA})$ and
$(\cB, L_{\cB})$ are isometric
if and only
if $\dist_{\cq}(\cA, \cB)=0$.
Then all the other assertions follow from
\cite[Theorem 1.1]{Li10}, Lemma~\ref{triangle dist_cq:lemma},
and Propositions~\ref{dist_cq<dist^R_cq:prop},
\ref{dist^R_oq<dist^R_cq:prop}, and \ref{dist^R_cq<dist^r_cq:prop}.

If $(\cA, L_{\cA})$ and $(\cB,
L_{\cB})$ are isometric, then clearly $\dist_{\cq}(\cA,\cB)=0$.
Suppose now $\dist_{\cq}(\cA, \cB)=0$.
By (\ref{dist_cq 1:eq}) we have $r_{\cA}=r_{\cB}$.
The case $r_{\cA}=0$ is trivial. So we assume that $r_{\cA}>0$.
We may also assume that both $(\cA, L_{\cA})$ and
$(\cB,L_{\cB})$ are closed. Set $X=\cD(A)$ and $Y=\cD(B)$.
For each $n\in \Ne$ we can find a complex normed space $V_n$ containing
$\cB$ and a linear isometric embedding
$h_n:\cA\hookrightarrow V_n$ such that $\dist^{V^{(3)}_n}_{\rH}(h^{(3)}_n(X^m),
Y^m)< \frac{1}{n}$. By \cite[Lemma 4.5]{Li10} we can find a complex
normed space $V$ containing all these $V_n$'s with the copies of
$\cB$ identified.
Then for any $n$ the union $\cup^{\infty}_{k=n+1}h^{(3)}_k(X^m)$ is
contained in the open $\frac{1}{n}$-neighborhood of $Y^m$, and hence
$\cup^{\infty}_{k=1}h^{(3)}_k(X^m)$ is
contained in the open $\frac{1}{n}$-neighborhood of $Y^m\cup(\cup^n_{k=1}h^{(3)}_k(X^m))$.
It follows immediately that $\cup^{\infty}_{k=1}h^{(3)}_k(X^m)$ is totally
bounded. Notice that $\{h^{(3)}_n|_{X^m}\}_{n\in \Ne}$ is an equicontinuous
sequence of mappings from $X^m$ to $V^{(3)}$. By the
Arzela-Ascoli theorem \cite{Conway90} there is a subsequence
$\{h^{(3)}_{n_k}|_{X^m}\}_{n\in \Ne}$ converging uniformly to a continuous map
$\varphi:X^m\to V^{(3)}$. By abuse of notation, we may assume
that this subsequence is $\{h^{(3)}_n|_{X^m}\}_{n\in \Ne}$ itself.
In particular, the sequence $\{h_n|_X\}_{n\in \Ne}$ converges
uniformly.
Clearly for every $a\in A+iA$ the sequence $\{h_n(a)\}_{n\in \Ne}$ also
converges. Since $A+iA$ is dense in $\cA$ and each $h_n$ is
isometric, it is easy to see
that  $\{h_n(a)\}_{n\in \Ne}$ converges for every $a\in \cA$.
Let $h(a)$ be the limit.
Clearly $h:\cA\to V$ is a linear isometric embedding.
Since $h_n|_X$ converges uniformly to $h|_X$, we have
$\dist^V_{\rH}(h_n(X), h(X))\to 0$. But
$\dist^V_{\rH}(h_n(X), Y)\le \dist^{V^{(3)}_n}_{\rH}(h^{(3)}_n(X^m),
Y^m)\to 0$.
So $\dist^V_{\rH}(Y, h(X))=0$, and hence
$Y=h(X)$.
Consequently $h(\cA)=\cB$. Clearly
$\varphi=h^{(3)}|_{X^m}$. Thus $h(aa')=h(a)h(a')$ for any $a, a'\in
X$. Consequently $h$ is an algebra homomorphism from $\cA$ to
$\cB$. Since $h(A)\subseteq B$, we have $h(\cA_{\sa})\subseteq
\cB_{\sa}$. Therefore $h:\cA\to \cB$ is a
$*$-isomorphism. In particular,
$h(e_{\cA})=e_{\cB}$. Then
\cite[Lemma 4.1]{Li10} tells us that $h(\cD_{\infty}(A))=\cD_{\infty}(B)$.
Hence $L_{\cA}=L_{\cA}\circ h$ on $\cA_{\sa}$.
So  $h$ is an isometry from $(\cA, L_{\cA})$ to
 $(\cB, L_{\cB})$.
\end{proof}

Next we show that our $C^*$-algebraic quantum distance
$\dist_{\cq}$ extends the ordinary Gromov-Hausdorff distance
$\dist_{\GH}$ in the way Rieffel's quantum distance $\dist_{\q}$ does (see
the discussion at the end of Section~\ref{Prelim:sec}).

\begin{theorem} \label{C*GH GH:thm}
The map $(X, \rho)\mapsto (C_{\Ce}(X),  L_{\rho})$ is a homeomorphism
  from
$(\CM,$ $\dist_{\GH})$ onto a closed subspace of
$(\CCQM,  \dist_{\cq})$.
\end{theorem}

\begin{proposition} \label{dist_cq<dist_GH:prop}
Let $(X, \rho_X)$ and $(Y, \rho_Y)$ be compact metric spaces.
For any $R\ge 0$ we have
\begin{eqnarray} \label{dist_cq<dist_GH:eq}
\dist^R_{\cq}(C_{\Ce}(X), C_{\Ce}(Y))\le \dist_{\GH}(X,
Y)\max(1, 2R).
\end{eqnarray}
\end{proposition}
\begin{proof}
The proof is similar to that of
\cite[Proposition 4.7 (8)]{Li10}.
Let $\rho$ be an admissible metric on $X\coprod Y$. Denote
$\dist^{\rho}_{\rH}(X, Y)$ by $d$. Denote by $Z$   the
set of elements $(x, y)$ in $X\times Y$ with
$\rho(x, y)\le d$. Since $X$ and $Y$ are compact, the
projections $Z\to X$ and $Z\to Y$ are surjective.
Thus the induced $*$-homomorphisms $h_X: C_{\Ce}(X)\to C_{\Ce}(Z)$ and
$h_Y: C_{\Ce}(Y)\to C_{\Ce}(Z)$ are faithful.
Notice that for any $f\in C_{\Ce}(X)$ and $g\in C_{\Ce}(Y)$ we have
\begin{eqnarray*}
\pa h_X(f)-h_Y(g)\pa=\sup\{|f(x)-g(y)|:(x, y)\in Z\}\le L_{\rho}(f, g)d,
\end{eqnarray*}
where $L_{\rho}$ is the Lipschitz seminorm on $C_{\Ce}(X\coprod Y)$.
Let $f_1, f_2\in \cD_R(C_{\Re}(X))$. Then we can extend $f_j$ to an
$F_j\in C_{\Re}(X\coprod Y)=C_{\Re}(X)\oplus C_{\Re}(Y)$ such that
$L_{\rho}(F_j)=L_{\rho_X}(f_j)$ and $\pa F_j\pa_{\infty}=\pa f_j\pa$,
where $\pa \cdot \pa_{\infty}$ is the supremum-norm on
$C_{\Ce}(X\coprod Y)$. For example, let
$F_j(w)=\min(\pa f_j\pa,\, \inf_{x\in X}(f_j(x)+L_{\rho_X}(f_j)\rho(x, w)))$ for all $w\in
X\coprod Y$. Say $F_j=(f_j, g_j)$ with $g_j\in C_{\Re}(Y)$.
Then $L_{\rho_Y}(g_j)\le L_{\rho}(F_j)\le 1$ and
$\pa g_j\pa \le  \pa F_j\pa_{\infty}=\pa f_j\pa\le R$.
So $g_j\in \cD_R(C_{\Re}(Y))$.
Then
\begin{eqnarray*}
\pa h_X(f_j)-h_Y(g_j)\pa &\le & L_{\rho}(F_j)d\le d, \\
L_{\rho}(f_1f_2, g_1g_2)&\le &
L_{\rho}(F_1)\pa F_2\pa_{\infty}+\pa F_1\pa_{\infty}L_{\rho}(F_2)
\le  2R,\\
\pa h_X(f_1f_2)-h_Y(g_1g_2)\pa &\le & L_{\rho}(f_1f_2, g_1g_2)d\le 2dR.
\end{eqnarray*}
Similarly, we can deal with pairs in
$\cD_R(C_{\Re}(X))$.
So $\dist^R_{\cq}(C_{\Ce}(X), C_{\Ce}(Y))\le \max(d,2dR)$ $=d\max(1, 2R)$.
Letting $\rho$ run over the admissible metrics on $X\coprod Y$,
we get (\ref{dist_cq<dist_GH:eq}).
\end{proof}

\begin{proof}[Proof of Theorem~\ref{C*GH GH:thm}]
The inequality (\ref{dist_cq<dist_GH:eq}) and Theorem~\ref{basic C*QGH:thm}
 tell us that this map is continuous.
Composing with the forgetful map $(\CCQM, \dist_{\cq})\rightarrow
(\CQM, \dist_{\oq})$ in Theorem~\ref{basic C*QGH:thm} we get also
a map $(\CM,  \dist_{\GH})\rightarrow (\CQM, \dist_{\oq})$.
Denote $(\CM,$
$\dist_{\GH})$ by $W_1$,  and denote its images in
$(\CCQM,  \dist_{\cq})$ and $(\CQM,  \dist_{\oq})$ by $W_2$ and $W_3$
respectively. Also denote the closure of $W_2$ in
$(\CCQM,  \dist_{\cq})$ by $\overline{W_2}$.
By \cite[Theorem 1.1]{Li10} $\dist_{\q}$ and $\dist_{\oq}$ are Lipschitz equivalent.
So we may identify $(\CQM,  \dist_{\oq})$ and $(\CQM,  \dist_{\q})$.
As pointed out
at the end of Section~\ref{Prelim:sec}, $W_3$ is
 closed in  $(\CQM,  \dist_{\oq})$, and the map $W_1\to W_3$ is a
 homeomorphism. Hence under the forgetful map
$(\CCQM,  \dist_{\cq})\to (\CQM,  \dist_{\oq})$
 the image of $\overline{W_2}$ is contained in $W_3$. Now
the composition $\overline{W_2}\to W_3\to W_1\to W_2$, where $W_3\to
 W_1$ is the inverse of $W_1\to W_3$,  is continuous and restricts to
 the identity map on $W_2$. Thus $\overline{W_2}=W_2$, and the map
 $W_1\to W_2$ is a homeomorphism with inverse $W_2\to W_3\to W_1$.
\end{proof}


\section{$C^*$-algebraic quantum Gromov compactness and completeness theorems}
\label{C*QGCC:sec}

In this section we prove the $C^*$-algebraic quantum Gromov
completeness and compactness theorems. Unlike in
Theorems~\ref{GH:thm} and Rieffel's quantum completeness and
compactness theorems \cite[Theorems 12.11 and 13.5]{Rieffel00} , we
have to put on some restriction. Just think about the
completeness. Take a Cauchy sequence $(\cA_n, L_n)$ in
$(\CCQM^R, \dist^R_{\cq})$. By (\ref{dist_cq 4:eq}) the
sequence $(A_n, L_n)$ is also Cauchy in $(\CQM^R,
\dist^R_{\oq})$. By \cite[Theorem 1.1]{Li10} and Rieffel's quantum completeness theorem
the space $(\CQM^R,
\dist^R_{\oq})$ is complete.
Denote by $(A, L)$ the limit of the
sequence $(A_n, L_n)$ in $(\CQM^R,
\dist^R_{\oq})$. Then $\cA:=\bar{A}+i\bar{A}$ should be the
limit of $(\cA_n, L_n)$, where $\bar{A}$ is the completion of $A$.
The trouble is how to get a
$C^*$-algebraic structure on $\cA$. For simplicity suppose
that $\cA_n$ and $\cA$ are all contained in a
normed space $V$ such that $\cD_R(A_n)$ converges to
$\cD_R(A)$ under $\dist^V_{\rH}$. Then for any $a, b\in \cD_R(A)$
we can find $a_n, b_n \in \cD_R(A_n)$ such that
$a_n \to a$ and $b_n \to b$ as $n\to \infty$. Clearly $a\cdot b$
should be defined as $\lim_{n\to \infty}a_nb_n$, if that limit
exists. Notice that $V$ is usually infinite-dimensional. The only way
we can make sure that $\{a_nb_n\}_{n\in \Ne}$ (or some subsequence of
$\{a_nb_n\}_{n\in \Ne}$) converges is that there exists a $\lambda \in
\Ce$ (not depending on $n$) such that $a_nb_n\in \lambda
\cD_R(A_n)$ for all $n$. Since $\pa a_nb_n\pa \le \pa
a_n\pa \cdot \pa b_n\pa \le R^2$, this is equivalent to say that $L_n(a_nb_n)$ is
uniformly bounded. In other words, $L_n(a_nb_n)$ should be controlled
by $L_n(a_n), L_n(b_n), \pa a_n\pa$, and $\pa b_n\pa$. We use Kerr's
definition \cite[page 155]{Kerr02}:

\begin{definition} \label{F Leib:def}
Let $F:\Re^4_+\to \Re_+$ be a continuous nondecreasing function, where $\Re^4_+$
is given the partial order $(x_1, x_2, x_3, x_4)\le (y_1, y_2, y_3,
y_4)$ if and only if $x_j\le y_j$ for all $j$. We say that a
$C^*$-algebraic compact quantum
metric space $(\cA, L)$ satisfies the \emph{$F$-Leibniz property}
if
\begin{eqnarray*}
L(a\cdot b)\le F(L(a), L(b), \pa a\pa, \pa b\pa)
\end{eqnarray*}
for all $a, b\in \cA$.
\end{definition}

\begin{remark} \label{F leib:remark}
(1) Rescaling $a$ and $b$ such that $\max(\pa a\pa, \, L(a))=\max(\pa
    b\pa, L(b))=1$, we see that any $(\cA, L)$ satisfying the
    $F$-Leibniz property for some $F$ also satisfies
     that $L(a\cdot b)\le
    F(1,1,1,1)(\pa a\pa+L(a))(\pa b\pa +L(b))$ for all $a, b\in \cA$;

(2)
The function $F(x_1, x_2, x_3, x_4)=x_1x_4+x_2x_3$ corresponds to the
Leibniz rule
\begin{eqnarray} \label{Leibniz:eq}
L(a\cdot b)&\le &L(a)\pa b\pa+\pa a\pa L(b).
\end{eqnarray}
The basic examples of $C^*$-algebraic compact quantum metric spaces
such as
$\theta$-deformations \cite{Li9}, quantum metric spaces induced by
ergodic actions \cite{Rieffel98b}, and
quantum metric spaces associated to ordinary
compact metric spaces
(see the discussion at the end of Section~\ref{Prelim:sec})
all satisfy (\ref{Leibniz:eq}).
But examples failing (\ref{Leibniz:eq})
also rise naturally.
A nonempty closed subset of a compact metric space is still such a space.
The quantum analogue also holds:  a nonzero quotient $\cB$
of a $C^*$-algebraic quantum
compact metric space $(\cA, L)$ is still such a space
(see \cite[Proposition 3.1]{Rieffel00} for the
corresponding assertion for compact quantum metric spaces).
But as the next lemma indicates,  even when $(\cA, L)$
satisfies (\ref{Leibniz:eq}),
we do not know whether $\cB$ satisfies it
or not.  See \cite{BC91} for similar results.
\end{remark}

\begin{lemma} \label{quotient of C*QCM:lemma}
Let  $(\cA, L_{\cA})$ be a
 $C^*$-algebraic compact quantum
metric space. For any unital
 $C^*$-algebra
$\cB$ and surjective unital positive map $\pi:\cA\to
 \cB$  let
$L_{\cB}$ be the quotient seminorm on $\cB$:
\begin{eqnarray*}
L_{\cB}(b):=\inf\{L_{\cA}(a):\pi(a)=b\}.
\end{eqnarray*}
 Then
$(\cB, L_{\cB})$ is a $C^*$-algebraic compact quantum
metric space. If $(\cA, L_{\cA})$ is closed, so is $(\cB,
 L_{\cB})$. Denote by $r$ the radius of $(\cA,
 L_{\cA})$.
If $\pi$ is a $*$-homomorphism (\ie $\cB$ is a quotient of
$\cA$) and $(\cA, L_{\cA})$ satisfies the
 $F$-Leibniz property, then $(\cB, L_{\cB})$ satisfies
 the $F_r$-Leibniz property, where $F_r(x_1, x_2, x_3, x_4)=F(x_1,
 x_2, 2x_3+4rx_1, 2x_4+4rx_2)$.
\end{lemma}
\begin{proof}
Clearly $L_{\cB}$ satisfies the reality condition
(\ref{real:eq}), and $L_{\cB}|_{\cB_{\sa}}$ is the
quotient of $L_{\cA}|_{\cA_{\sa}}$.
Since the subspace of $\cA_{\sa}$
with finite $L_{\cA}$ is a compact quantum metric space
and is dense in $\cA_{\sa}$,
by \cite[Proposition 3.1]{Rieffel00} the subspace of
$\cB_{\sa}$ with finite $L_{\cB}$ is a compact quantum
metric space and is dense in $\cB_{\sa}$.
Thus $(\cB, L_{\cB})$ is a $C^*$-algebraic compact quantum
metric space.
When $(\cA, L_{\cA})$ is closed,
an argument similar to
that in the proof of \cite[Proposition 3.3]{Rieffel00} shows that
$(\cB, L_{\cB})$ is also closed.

Now suppose that $\pi$ is a $*$-homomorphism and that
$(\cA, L_{\cA})$ satisfies the
 $F$-Leibniz property.
For any $b_1, b_2\in \cB$ and $\varepsilon>0$, take $a_1, a_2\in
\cA$ such that
$\pi(a_j)=b_j$ and $L_{\cA}(a_j)\le
L_{\cB}(b_j)+\varepsilon$. Denote by $(a_j)_{\sa}$ and
$(a_j)_{\im}$ the self-adjoint part and the imaginary part of $a_j$ respectively.
From (\ref{real:eq}) we have $L_{\cA}((a_j)_{\sa}), \, L_{\cA}((a_j)_{\im}) \le L_{\cA}(a_j)$.
By \cite[Lemma 3.4]{Rieffel00} we get
$\pa (a_j)_{\sa}\pa,\,  \pa (a_j)_{\im}\pa
\le \pa b_j\pa +
2rL_{\cA}(a_j)$. So
$\pa a_j\pa \le
2\pa b_j\pa + 4r(L_{\cB}(b_j)+\varepsilon)$. Notice that
$\pi(a_1a_2)=b_1b_2$. Therefore
\begin{eqnarray*}
& & L_{\cB}(b_1b_2) \\
&\le & L_{\cA}(a_1a_2)\le
F(L_{\cA}(a_1),L_{\cA}(a_2), \pa a_1\pa, \pa a_2\pa) \\
&\le &
F(L_{\cB}(b_1)+\varepsilon, L_{\cB}(b_2)+\varepsilon, 2\pa
b_1\pa + 4r(L_{\cB}(b_1)+\varepsilon), 2\pa b_2\pa + 4r(L_{\cB}(b_2)+\varepsilon)).
\end{eqnarray*}
Letting $\varepsilon\to 0$, we get
\begin{eqnarray*}
L_{\cB}(b_1b_2)
\le F_r(L_{\cB}(b_1), L_{\cB}(b_2), \pa b_1\pa, \pa b_2\pa).
\end{eqnarray*}
\end{proof}

Clearly if $(\cA, L)$ satisfies the F-Leibniz property, then so
does its closure. Denote by $\CCQM_F$
the set of all isometry classes of
$C^*$-algebraic compact quantum metric spaces
satisfying the F-Leibniz
property.

\begin{theorem}[$C^*$-algebraic Quantum Gromov Completeness and
  Compactness Theorems] \label{C*QGH:thm}
Let $F:\Re^4_+\to \Re_+$ be a continuous nondecreasing function.
Then $(\CCQM_F, \dist_{\cq})$ is a complete
metric space. A subset
$\mathcal{S}\subseteq \CCQM_F$ is totally bounded if and only if

(1) there is a constant $D$ such that $\diam(\cA, L)\le D$ for all
$(\cA, L)\in \mathcal{S}$;

(2) For any $\varepsilon>0$, there exists a constant $K_{\varepsilon}>0$ such that
$\Cov(\cA, \varepsilon)\le K_{\varepsilon}$ for all
$(\cA, L)\in \mathcal{S}$,

\noindent
and if and only if

(1') there is a constant $D'$ such that
$\diam(\cD_{1,r_{\cA}}(A))\le D'$ for all
$(\cA, L)\in \mathcal{S}$;

(2') For any $\varepsilon>0$, there exists a constant $K'_{\varepsilon}>0$
such that $\Cov(\cD(A), \varepsilon)\le K'_{\varepsilon}$ for all $(\cA, L)\in
\mathcal{S}$.

\end{theorem}

\begin{notation} \label{X_R:notation}
Let $R\ge 0$. For any compact metric space $(X, \rho)$ let
$C(X)_R:=\{f\in C_{\Ce}(X):L_{\rho}(X)\le 1, \pa f\pa \le R\}$, equipped with the
metric induced from the supremum-norm in $C_{\Ce}(X)$, where $L_{\rho}$
is the Lipschitz seminorm defined by (\ref{dist to Lip:eq}).
\end{notation}

\begin{lemma}
 \label{univ embedQ:lemma}
Let $\mathcal{S}$ be a subset of $\CCQM$
satisfying the conditions (1) and (2) in Theorem~\ref{C*QGH:thm}.
Let $R\ge \sup\{r_{\cA}:(\cA, L)\in \mathcal{S}\}$.
Then there exist a complex Banach space $V$ and
a compact convex subset $Z\subseteq V$ such that for
every $(\cA, L)\in \mathcal{S}$ there is
a linear isometric embedding $h_{\cA}:
\cA\hookrightarrow V$ with
 $h_{\cA}(\cD_R(A))\subseteq Z$.
\end{lemma}
\begin{proof}
We may assume that every $(\cA, L)$ is closed, and that
$R>0$. By \cite[Lemma 5.4]{Li10}
the set $\{\cD_R(A):(\cA, L)\in \mathcal{S}\}$
satisfies the condition (2) in Theorem~\ref{GH:thm}.
Let $(\cA, L)\in \mathcal{S}$. Denote by
$\De(\cA')$  the unit ball of
$\cA'$.
Then  $\Ce\cD_{R}(A)=A+iA$ is dense in
  $\cA$, and hence  the natural map
  $\psi_{\cA}:\De(\cA')\to
  C_{\Ce}(\cD_{R}(A))$ induced by the pairing
  between $\cA$ and $\cA'$ is injective. Clearly
  $\psi_{\cA}$ maps $\De(\cA')$ into
  $C(\cD_R(A))_R$
(see Notation~\ref{X_R:notation}).
Endow
$\De(\cA')$ with the metric $\rho$ induced by $\psi_{\cA}$ and the
supremum-norm in $C_{\Ce}(\cD_R(A))$. Then
\cite[Corollary 5.3, Lemma 5.1]{Li10} tell
us that the set $\{\De(\cA'): (\cA, L)\in
\mathcal{S}\}$ satisfies the condition (2) in Theorem~\ref{GH:thm}.
By \cite[Proposition 5.2]{Li10} we can find a complex Banach space $V$ and
a compact convex subset $Z\subseteq V$ such that for
every $(\cA, L)\in \mathcal{S}$ there is
a linear isometric embedding $\phi_{\cA}:C_{\Ce}(\De(\cA'))
\hookrightarrow V$ with the image of $C(\De(\cA'))_R$
contained in $Z$.
Notice that $\rho(f, g)=\sup\{|f(a)-g(a)|:a\in
\cD_R(A)\}$ for any $f, g\in\De(\cA')$.
It is easy to see that the topology defined
by $\rho$ on $\De(\cA')$ is exactly the
weak-$*$ topology.
Then the natural pairing between $\cA$ and $\cA'$
also gives a linear isometric embedding $\varphi_{\cA}:\cA\hookrightarrow
C_{\Ce}(\De(\cA'))$. Clearly $\varphi$ maps
$\cD_R(A)$
into $C(\De(\cA'))_R$.
Thus we may just set $h_{\cA}=\phi_{\cA}\circ \varphi_{\cA}$.
\end{proof}

\begin{proof}[Proof of Theorem~\ref{C*QGH:thm}]
We prove the compactness part first.
The equivalence between (1)+(2) and (1')+(2') follows from
\cite[Lemma 5.4]{Li10}.
Suppose that $\mathcal{S}$ is totally bounded. By (\ref{dist_cq 1:eq})
the sets $\{r_{\cA}:(\cA, L)\in \mathcal{S}\}$ and
$\{\cD(A):(\cA, L)\in
\mathcal{S}\}$ are totally bounded in $\Re$ and $(\CM,
\dist_{\GH})$ respectively. Then the condition (1) follows
immediately. Also the condition (2) follows from
Theorem~\ref{GH:thm}.

Now suppose that $\mathcal{S}$ satisfies
the conditions (1) and (2). We shall show that $\mathcal{S}$ is
pre-compact, which implies that $\mathcal{S}$ is totally bounded.
Set $R=1+\sup\{r_{\cA}:(\cA, L)\in \mathcal{S}\}$.
By Lemma~\ref{univ embedQ:lemma} we can find a complex Banach
space $V$ and a compact convex subset $Z_1\subseteq V$ such that $\cA$ is a linear
subspace of $V$ with $\cD_R(A)\subseteq Z_1$ for all
$(\cA, L)\in \mathcal{S}$. Set $R'=\max(R, F(1,1,R,R))$, and
set $Z=R'Z_1+iR'Z_1$. Then $Z\supseteq Z_1$. Let $(\cA, L)\in
\mathcal{S}$, and let $a, b\in \cD_R(A)$.
Then $L(ab)\le F(L(a), L(b), \pa a\pa, \pa
b\pa)\le F(1, 1, R, R)\le R'$. Also $\pa ab\pa \le R^2\le R\cdot
R'$. Thus the self-adjoint and imaginary parts of $ab$ are both
in $R'\cD_R(A)\subseteq R'Z_1$. Then $ab\in Z$.
 Therefore
$(\cD_R(A))^m$ is contained in $Z^{(3)}:=Z\oplus Z\oplus
Z\subseteq V^{(3)}$. Notice that $Z^{(3)}$ is compact.
Then the set of nonempty closed subsets of $Z^{(3)}$ is compact equipped with
the Hausdorff distance $\dist^{Z^{(3)}}_{\rH}$.
Thus we can find a sequence
$\{(\cA_n, L_n)\}_{n\in \Ne}$ in $\mathcal{S}$ such that
$(\cD_R(A_n))^m$ converges under $\dist^{V(3)}_{\rH}$. Denote by $Y$
the limit. Set $\pi_j:V^{(3)}\to V$ to be the projection of $V^{(3)}$
to the $j$-th coordinate. Also set $\pi_{12}:V^{(3)}\to V^{(2)}$ to be the projection of $V^{(3)}$
to the first two coordinates. Clearly $\dist^V_{\rH}(\cD_{R}(A_n), \pi_j(Y))
=\dist^{V}_{\rH}(\pi_j((\cD_R(A_n))^m),
  \pi_j(Y))\le \dist^{V^{(3)}}_{\rH}((\cD_R(A_n))^m,
  Y)$ for $j=1,2$. Thus $\cD_{R}(A_n)$ converges to $\pi_1(Y)$ and $\pi_2(Y)$ under
$\dist^V_{\rH}$. Set $X=\pi_1(Y)=\pi_2(Y)$. Similarly,
$(\cD_R(A_n))^{(2)}=\pi_{12}((\cD_{R}(A_n))^m)$ converges to $\pi_{12}(Y)$ under
$\dist^{V^{(2)}}_{\rH}$. But clearly $(\cD_R(A_n))^{(2)}$
converges to $X^{(2)}$ under $\dist^{V^{(2)}}_{\rH}$. Thus
$\pi_{12}(Y)=X^{(2)}$.

Since each $\cD_R(A_n)$ is $\Re$-balanced
(\ie $\lambda a\in \cD_R(A_n)$ for all $a\in \cD_R(A_n)$ and $\lambda \in \Re$
with $|\lambda|\le 1$)
and convex, has radius $R$ and contains $0_V$,
clearly so does $X$. Thus the set $\Re_+\cdot X=\{\lambda x:\lambda
\in \Re_+, x\in X\}$ is a real linear subspace of $V$. Denote it by
$B$, and denote the closure of $B+iB$ by $\cB$. Then
$\cB$ is a closed complex linear subspace of $V$.
We shall define a $C^*$-algebra structure and a Lip-norm on $\cB$.
It is easy to see that for every $a\in \cB$ we can find
$a_n\in \cA_n$ for each $n$ such that $a_n\to a$ as $n\to \infty$.

\begin{lemma} \label{prod:lemma}
Let $a, b\in \cB$.
Let $a_n, b_n\in \cA_n$ for each $n$ such that $a_n\to a$ and
$b_n\to b$ as $n\to \infty$. Then the sequence $\{a_nb_n\}_{n\in \Ne}$
converges to an element in $\cB$,
and the limit depends only on $a$ and $b$.
\end{lemma}
\begin{proof}
We show first that for any $x,y \in B+iB$ there exist
$x_n, y_n\in \cA_n$ for each $n$ such that $x_n\to x, \, y_n\to y$, and
$\{x_ny_n\}_{n\in \Ne}$
converges to an element in $\cB$.
Since $B=\Re_+\cdot X$, it suffices to show this for all
$x, y\in X$.
Let $x, y\in X$. Since $\pi_{12}(Y)=X^{(2)}$, we can find $z\in V$
with $(x, y, z)\in Y$. Then we can pick $(x_n, y_n, x_ny_n)\in
(\cD_R(A_n))^m$ for each $n$ such that $(x_n, y_n, x_ny_n)\to (x, y,
z)$ in $V^{(3)}$ as $n\to \infty$. Notice that
$x_ny_n\in R'\cD_R(A_n)+iR'\cD_R(A_n)$.
Thus $\lim_{n\to \infty}(x_ny_n)$ is in $R'X+iR'X\subseteq
\cB$.

Now let $a$ and $b$ be arbitrary elements in $\cB$, and let
$a_n, b_n\in \cA_n$ with $a_n\to a$ and $b_n\to b$.
Let $\varepsilon>0$. Pick $x, y\in B+iB$ with $\pa a-x\pa, \pa b-y\pa <
\varepsilon$. Let $x_n$ and $y_n$ be as in the above.
Take $N_1$ such that $\pa a_n-a\pa, \pa b_n-b\pa, \pa x_n-x\pa, \pa
y_n-y\pa<\varepsilon$ for all $n>N_1$.
Since  $\pa
a_n-x_n\pa\to \pa a-x\pa$ and $\pa b_n-y_n\pa \to \pa b-y\pa$, there
exists $N_2>N_1$ such that $\pa a_n- x_n\pa, \pa b_n-y_n\pa
<\varepsilon$ for all $n>N_2$. Then
\begin{eqnarray*}
\pa a_nb_n-x_ny_n\pa &\le & \pa a_n-x_n\pa \cdot \pa
b_n\pa +\pa x_n \pa \cdot \pa b_n -y_n\pa\\
&\le &\varepsilon (\pa b\pa
+\varepsilon)+(\pa x\pa +\varepsilon)\varepsilon
\le \varepsilon(\pa b\pa +\pa
a\pa+3\varepsilon)
\end{eqnarray*}
for all $n>N_2$.
Since $\{x_ny_n\}_{n\in \Ne}$ converges, we can find $N_3>N_2$ such
that $\pa x_ny_n-x_ky_k\pa \le \varepsilon$ for all $n, k>N_3$.
Then $\pa a_nb_n-a_kb_k\pa \le \varepsilon+2\varepsilon(\pa b\pa +\pa
a\pa+3\varepsilon)$ for all $n, k>N_3$. Thus $\{a_nb_n\}_{n\in \Ne}$
is a Cauchy sequence, and hence converges. We also get that
$\pa \lim_{n\to \infty}(a_nb_n)-\lim_{n\to \infty}(x_ny_n)\pa \le
\varepsilon(\pa b\pa +\pa a\pa+3\varepsilon)$. Since $\lim_{n\to
  \infty}(x_ny_n)$ does not depend on the choice of the sequences
$\{a_n\}_{n\in \Ne}$ and
$\{b_n\}_{n\in \Ne}$, neither does $\lim_{n\to
  \infty}(a_nb_n)$. The above inequality also shows that $\lim_{n\to
  \infty}(a_nb_n)$ is in $\cB$.
\end{proof}
For any $a, b\in \cB$ denote by $a\cdot b$ the limit
$\lim_{n\to \infty}(a_nb_n)$ in Lemma~\ref{prod:lemma}. Then clearly
this makes
$\cB$ into an algebra over $\Ce$. Also notice that
\begin{eqnarray*}
\pa a\cdot
b\pa =\pa \lim_{n\to \infty}(a_nb_n)\pa \le \lim_{n\to \infty}\pa a_n
\pa \cdot \pa b_n\pa=\pa a\pa \cdot \pa b\pa.
\end{eqnarray*}
Thus $\cB$ is a
Banach algebra. It also follows from the proof of
Lemma~\ref{prod:lemma}
that $Y=X^m$. By Lemma~\ref{cq o-u:lemma} the sequence
$\{e_{\cA_n}\}_{n\in \Ne}$ converges. Denote by $e_{\cB}$
the limit. Clearly $e_{\cB}$ is an identity of $\cB$.

For any $a_n\in \cA_n$ with $a_n\to a\in
\cB$, by a similar argument as in Lemma~\ref{prod:lemma} we
see that $\{a^*_n\}_{n\in \Ne}$ converges to an element in
$\cB$. Denote the limit by $a^*$. Then this makes
$\cB$ into a $*$-algebra. Clearly
$\bar{B}=(\cB)_{\sa}$.
Notice that
\begin{eqnarray*}
\pa a^*a\pa=\pa \lim_{n\to \infty} a^*_na_n\pa =\lim_{n\to \infty}
\pa a^*_na_n\pa=\lim_{n\to \infty}
\pa a_n\pa^2=\pa a\pa^2.
\end{eqnarray*}
Hence $\cB$ is $C^*$-algebra.

We proceed to define a Lip-norm $L$ on $\cB$.
As hinted by \cite[Lemma 4.1]{Li10}, we describe $L|_B$
first. For this we need:
\begin{lemma} \label{X+R:lemma}
We have $X=\{b\in (X+\Re e_{\cB}):\pa b\pa \le R\}$.
\end{lemma}
\begin{proof}
Clearly we have $X\subseteq \{b\in (X+\Re e_{\cB}):\pa b\pa \le R\}$.
Let $x\in X$ and $\lambda \in \Re$ with  $\pa x+\lambda e_{\cB}\pa \le
R$. Pick $x_n\in \cD_R(A_n)$ for each $n\in \Ne$
such that $x_n\to x$. Then $x_n+\lambda e_{\cA_n}\to x+\lambda
e_{\cB}$, and hence $\lim_{n\to \infty}\pa x_n+\lambda
e_{\cA_n}\pa \le R$. Set $t_n=\max(R, \pa x_n+\lambda
e_{\cA_n}\pa)$, and set $y_n=\frac{R}{t_n}( x_n+\lambda
e_{\cA_n})$. Then $y_n\in \cD_R(A_n)$, and
$y_n\to x+\lambda
e_{\cB}$. Thus $x+\lambda
e_{\cB}\in X$.
\end{proof}
By \cite[Lemma 4.1]{Li10} there is a unique closed Lip-norm
$L_*$ on $B$ satisfying $\cD_R(B)=X$. And $(B, L_*)$ has
radius no bigger than $R$.

Now we describe the whole $L$. As hinted by
\cite[Proposition 6.9]{Li10}, we have:

\begin{lemma} \label{lip:lemma}
Define a (possibly $+\infty$-valued) seminorm, $L$, on $\cB$ by
\begin{eqnarray} \label{limit L:eq}
L(a):=\inf\{\limsup_{n\to \infty}L_n(a_n): a_n\in \cA_n \mbox{
  for all } n, \mbox{ and } a_n\to a\}.
\end{eqnarray}
Then $(\cB, L)$ is $C^*$-algebraic compact quantum metric space
with radius no bigger than $R$ and satisfying the $F$-Leibniz
property.
Also $L|_B=L_*$.
\end{lemma}
\begin{proof}
Clearly $L$ is a seminorm satisfying the reality condition
(\ref{real:eq}). Also clearly $L(x)\le 1$ for every $x\in X$, and
$L(e_{\cB})=0$. Thus $\Re e_{\cB}+X\subseteq \{a\in \cB: L(a)\le
1\}$. It is easy to see that the infimum in the definition of $L$
can always be achieved. Let $a\in \cB_{\sa}$ with $L(a)\le 1$.
Then there exist $a_n\in  \cA_n$ such that $a_n\to a$ and
$\limsup_{n\to \infty}L_n(a_n)=L(a)\le 1$. Replacing $a_n$ by the
self-adjoint part of $a_n$, we may assume that $a_n\in
(\cA_n)_{\sa}$. Rescaling $a_n$ slightly, we may also assume that
$L_n(a_n)\le L(a)$ and $\pa a_n\pa \le \pa a\pa$ for all $n\in
\Ne$. By \cite[Lemma 4.1]{Li10} we can write $a_n$ as
$x_n+\lambda_n e_{\cA_n}$ for some $x_n\in \cD_R(A_n)$ and
$\lambda_n\in \Re$. Then $|\lambda_n|\le \pa a_n\pa +\pa x_n\pa
\le \pa a\pa +R$. Thus some subsequence $\{\lambda_{n_k}\}_{k\in
\Ne}$ converges to a $\lambda \in \Re$. Since $\{  x_n+\lambda_n
e_{\cA_n}\}_{n\in \Ne}$ converges, the subsequence
$\{x_{n_k}\}_{k\in \Ne}$ also converges. Denote by  $x$ the limit.
Then $x\in X$ and $a=x+\lambda e_{\cB}$. Therefore $\Re e_{\cB}+X=
\{a\in \cB_{\sa}: L(a)\le 1\}$. By \cite[Lemma 4.1]{Li10} we have
$L|_B=L_*$. Hence $(\cB, L)$ is $C^*$-algebraic compact quantum
metric space with radius no bigger than $R$. Since each $(\cA_n,
L_n)$ satisfies the F-Leibniz property, clearly so does $(\cB,
L)$.
\end{proof}

Back to the proof of Theorem~\ref{C*QGH:thm}.
We have seen that $(\cD_R(A_n))^m$ converges to $Y=X^m=(\cD_R(B))^m$
under $\dist^{V^{(3)}}_{\rH}$. Thus
$$\dist^R_{\cq}(\cA_n,
\cB)\le \dist^{V^{(3)}}_{\rH}((\cD_R(A_n))^m, X^m)\to
0.$$
By Theorem~\ref{basic C*QGH:thm} we get $\dist_{\cq}(\cA_n,
\cB)\to 0$. This shows that $\mathcal{S}$ is pre-compact.

We are left to show that $(\CCQM_F, \dist_{\cq})$ is
complete.   Let  $\{(\cA_n, L_n)\}_{n\in \Ne}$
be a Cauchy sequence in $(\CCQM_F, \dist_{\cq})$. From
(\ref{dist_cq 1:eq}) and Theorem~\ref{GH:thm}
one sees that $\mathcal{S}:=\{(\cA_n, L_n):n\in
\Ne\}$ satisfies the conditions (1') and (2'). Thus $\mathcal{S}$ is
pre-compact. Then some subsequence of $\{(\cA_n, L_n)\}_{n\in
  \Ne}$
converges under $\dist_{\cq}$. Since $\{(\cA_n, L_n)\}_{n\in
  \Ne}$ is Cauchy, it converges itself. This finishes the proof
of Theorem~\ref{C*QGH:thm}.
\end{proof}

\begin{remark} \label{C*QGH:remark}
Similarly, one can give a direct proof of \cite[Theorem 5.5]{Li10}.
One can also give a simple proof of Rieffel's
quantum completeness and compactness theorems (as formulated
in \cite[Theorem 2.4]{Li10}) along the
same line: the proof for the completeness and the ``only if'' part is
essentially the same as in the proof of Theorem~\ref{C*QGH:thm}. We
indicate briefly how to show that a subset $\mathcal{S}\subseteq
\CQM$ satisfying the conditions (1) and (2) is precompact:
Pick $R>\sup\{r_A:(A, L)\in \mathcal{S}\}$. Then the natural pairing
between $A$ and $A'$ gives an affine isometric embedding
$S(A)\hookrightarrow C(\cD_R(A))_R$ (see Notation~\ref{X_R:notation}).
By \cite[Lemma 5.4, Proposition 5.2]{Li10} we can find a normed space $V$
and a compact convex subset $Z$ such that $Z\supseteq
C(\cD_R(A))_R\supseteq S(A)$ for all $(A, L)\in
\mathcal{S}$.
Since the set of convex closed subsets of $Z$ is compact equipped with
the Hausdorff distance $\dist^Z_{\rH}$,
some sequence
$\{S(A_n)\}_{n\in \Ne}$ converges to a convex subset $X\subseteq Z$
under $\dist^V_{\rH}$. By \cite[Proposition 3.1]{Li10} $X$ is the
state-space of some compact quantum metric  space $(B, L)$. Then
\cite[Proposition 3.2]{Li10} tells us that $(A_n, L_n)\to (B, L)$
under $\dist_{\q}$.
\end{remark}

\begin{remark} \label{Weaver:remark}
In Definition~\ref{C*QGH:def} we defined our $C^*$-algebraic quantum
distance as
$\inf\{\dist^{V^{(3)}}_{\rH}((\cD(A))^m,
(\cD(B))^m)\}$, where $V$ runs over complex
normed spaces containing $\cA$ and $\cB$.
 From the reality condition
(\ref{real:eq}) one can see easily that the ball
$\cD(\cA)$ is also totally
bounded.  Then one can try to define
a quantum distance, $\dist'_{\cq}$,  as
$\inf\{\dist^{V^{(3)}}_{\rH}((\cD_{\diam(\cA)}(\cA))^m,
(\cD_{\diam(\cB)}(\cB))^m)\}$.
One can show that $\dist'_{\cq}$ is still
 a metric on the set of  isometry classes of
$C^*$-algebraic
compact quantum metric spaces, though now the isometry in
Definition~\ref{iso of C*QCM:def} should be required to preserve the
Lip-norm on the whole algebra $\cA$
instead of only on $\cA_{\sa}$. At first sight, this sounds
more satisfactory. However,  the map
$(X, \rho)\mapsto (C_{\Ce}(M),  L_{\rho})$ turns out not to be continuous with
respect to $\dist'_{\cq}$. This phenomenon is rooted in the different
extension
behavior of $\Ce$-valued Lipschitz functions and
$\Re$-valued Lipschitz functions on metric spaces:
for any metric space $(X, \rho)$ and a
subspace $Y\subseteq X$, if $f$ is an $\Re$-valued Lipschitz function
on $Y$, we can always extend $f$ to a function $\hat{f}$ (for
  example, $\hat{f}(x)=\inf_{y\in Y}(f(y)+L(f)\rho(x, y))$ or
$\hat{f}(x)=\sup_{y\in Y}(f(y)-L(f)\rho(x, y))$) on $X$ with
the same Lipschitz seminorm. But when $f$ is $\Ce$-valued, sometimes one
has to increase the Lipschitz seminorm to extend $f$. The following
example comes from \cite{Weaver99}:

\begin{example}\cite[Example 1.5.7]{Weaver99} \label{4 points:eg}
Let $X$ be a four-element set, say $X=\{x, y_1$, $y_2, y_3\}$, with $\rho(x,
y_i)=\frac{1}{2}$ for all $i$ and $\rho(y_j, y_k)=1$ for distinct
$j,k$. Let $Y=X\setminus \{x\}$, and let $f:Y\rightarrow \Ce$ be an
isometric map. Then $L(f)=1$ and every extension $\hat{f}$ of $f$ to
$X$ has Lipschitz seminorm at least $\frac{2}{\sqrt{3}}$.
\end{example}

The nice behavior of $\Re$-valued functions with respect to Lipschitz
seminorms is also reflected in
\cite[Proposition 6.11]{Li10}.
In fact, the analogue of \cite[Proposition 6.11]{Li10} for
$\Ce$-valued functions
holds for convergence under $\dist'_{\cq}$ (at least for commutative
case), but not for convergence under $\dist_{\GH}$.
 This explains the discontinuity of the map
$(X, \rho)\mapsto (C_{\Ce}(M),  L_{\rho})$ with
respect to $\dist'_{\cq}$.
Let us give one example.
Suppose that $\{X_n\}_{n\in \Ne}$ is a sequence of  closed subsets of a compact metric space
$Z$ and that $X_n$ converges to some $X_0\subseteq Z$ under
$\dist^Z_{\rH}$. Then the restrictions of elements in $C_{\Ce}(Z)$ generate
a continuous field of $C^*$-algebras over $T=\{\frac{1}{n}:n\in
\Ne\}\cup \{0\}$ with fibres $C_{\Ce}(X_n)$ at $t=\frac{1}{n}$ and
$C_{\Ce}(X_0)$ at $t=0$. Denote by  $\Gamma$ the space of continuous
sections. Set $L_n=L_{\rho_{X_n}}$, and set $L_0=L_{\rho_{X_0}}$.
As in the proof of Theorem~\ref{C*QGH:thm}
embed all the $C_{\Ce}(X_n)$'s isometrically into a common normed space
$V$.
Notice that the proof of
Theorem~\ref{C*QGH:thm} can be easily modified to construct the limit
$(\cB, L)$
of some subsequence of  $\{(C_{\Ce}(X_n), L_n)\}_{k\in \Ne}$ under
$\dist'_{\cq}$. Then $\cB$ is a subspace of $V$,
and $L$ is still given by (\ref{limit L:eq}) of Lemma~\ref{lip:lemma}.
Without loss of generality, we may assume that this subsequence
is $\{C_{\Ce}(X_n)\}_{n\in \Ne}$ itself.
Take a dense sequence $\{a_m\}_{m\in \Ne}$ in $\cD_{R}(C_{\Ce}(X_0))$,
and take $f_m\in \Gamma$ with $(f_m)_0=a_m$.
By a standard diagonal argument we can find a subsequence of
$\{C_{\Ce}(X_n)\}_{n\in \Ne}$ such that the corresponding
subsequence of $\{(f_m)_{1/n}\}_{n\in \Ne}$ converges in $V$ for every
$m$. Without loss of generality, we still assume that this subsequence
is $\{C_{\Ce}(X_n)\}_{n\in \Ne}$ itself. Notice that the limit
$b_m=\lim_{n\to \infty}(f_m)_{1/n}$ is in $\cB$, and that the
map $a_m\mapsto b_m$ extends to an isometric affine map from $\cD_{R}(C_{\Ce}(X_0))$
into $\cB$. Then one can verify that
this map extends to a $*$-isomorphism from $C_{\Ce}(X_0)$ onto $\cB$.
Identify $C_{\Ce}(X_0)$ with $\cB$. Then one verifies that
elements in $\Gamma$ are exactly continuous maps $\psi: T\to V$
with $\psi(\frac{1}{n})\in C_{\Ce}(X_n)$ and $\psi(0)\in
C_{\Ce}(X_0)$. Thus (\ref{limit L:eq}) becomes
$L(a)=\inf\{\limsup_{n\to \infty}L_n(f_{\frac{1}{n}}): f\in \Gamma
  \mbox{ and } f_0=a\}$, which is the analogue of
\cite[Proposition 6.11]{Li10}.
As pointed out in the proof of
Lemma~\ref{lip:lemma},
the infimum here can always be achieved
  by some $f\in \Gamma$.

I am indebted to Nik Weaver showing the next example \cite{Weaver01}.
Note that this example
modifies Example~\ref{4 points:eg}, and shows that
the limit Lip-norm $L$ so constructed does not coincide with
$L_0$ on the whole of $C_{\Ce}(X_0)$. Consequently, the map
$(X, \rho)\mapsto (C_{\Ce}(M),  L_{\rho})$ is not continuous with
respect to $\dist'_{\cq}$.

\begin{example} \label{modified 4 point:eg}
Let $X_0$ be a solid equilateral triangle in $\Ce$
with sides of length $1$ and center
at the origin. The metric $\rho_{X_0}$ is just the Euclidean metric on
$\Ce$.
 We define $X_n$
as follows: subdivide $X_0$ into $4^n$
congruent equilateral sub-triangles with sides of length $1/2^n$.
Let $V_n$
be the vertices of these sub-triangles, and let $B_n$ be their barycenters.
Then let $X_n$ be $V_n\cup B_n$. For each of the sub-triangles
put edges of length $1/2^n$ between its
vertices and put edges of length $1/2^{n+1}$ between each
vertex and the barycenter.  Now extend the distance on
$V_n$ to a distance on $X_n$ via setting distances between points in $B_n$ and
points in $X_n$ to be the
graph distance, \ie $d(x,y) =$ length of the shortest path
between $x$ and $y$ for any $x\in B_n$ and any $y\in X_n$.
We construct a metric space $Z$ containing all
these $X_n$'s and $X_0$ such that $X_n\to X_0$ under $\dist^Z_{\rH}$:

\begin{lemma} \label{glue:lemma}
There exists a compact metric space $(Z, \rho)$ with isometric
embeddings $\varphi_n:X_n\hookrightarrow Z$ for all $n\in \Ne$ and
$\varphi_0:X_0\hookrightarrow Z$ such that $\varphi_m(p)=\varphi_0(p)$
for each $p\in V_n$ and all $m\ge n$.
\end{lemma}
\begin{proof}
Notice that for any metric space $X$ and $Y$ with a common closed
subspace $W$, we can ``glue'' $X$ and $Y$ together along $W$ by
setting $\rho(x, y)=\inf_{w\in W}(\rho_X(x, w)+\rho_Y(w, y))$ for
all $x\in X,  y\in Y$.
Gluing $X_0$ and $X_1$ together along $V_1$, we get a metric on
the disjoint union $X_0\coprod B_1$. Gluing $X_0\coprod B_1$ and $X_2$ together along $V_2$
we get a metric on $X_0\coprod B_1\coprod B_2$. In this way, we get a
metric on $X_0\coprod(\coprod_{n\in N}B_n)$. Clearly
$X_0\coprod(\coprod_{n\in N}B_n)$ is totally bounded. So the
completion of $X_0\coprod(\coprod_{n\in N}B_n)$, denoted by $Z$, is
compact. The natural embeddings $\varphi_n:X_n\hookrightarrow Z$
and $\varphi_0:X_0\hookrightarrow Z$ satisfy our requirement.
\end{proof}

Identify $X_n$ and $X_0$ with their images in $Z$. Then
\begin{eqnarray*}
\dist^Z_{\rH}(X_0, X_n)&\le &\dist^Z_{\rH}(X_0, V_n)+\dist^Z_{\rH}(V_n, X_n)\\&\le &
\dist^{X_0}_{\rH}(X_0, V_n)+\dist^{X_n}_{\rH}(V_n, X_n)\to 0.
\end{eqnarray*}
Let $a$ be
the identity map from $X_0$ into $\Ce$.
Then we have:

\begin{proposition} \label{modified 4 point:prop}
For any $a_n:X_n\rightarrow \Ce$ with  $L_n(a_n)\le 1$,
we have $\pa a_n|_{V_n} - a|_{V_n}\pa  \geq \frac{\sqrt{3}}{6} - \frac{1}{4}$.
\end{proposition}
\begin{proof}
Suppose that $\pa a_n|_{V_n} - a|_{V_n}\pa < \frac{\sqrt{3}}{6} - \frac{1}{4}$.
In particular, for any $x$ in $V_n$ that lies on the boundary of $X_0$,
we have $|a_n(x) - a(x)| <\frac{\sqrt{3}}{6} - \frac{1}{4}$.  It follows by some
elementary geometry that $\frac{\sqrt{3}}{2}X_0$ is strictly contained in the
closed curve formed by joining the points $a_n(x)$ (for $x$ in
$V_n\cap \partial X_0$, where $\partial X_0$ is the boundary of
$X_0$) with line segments.

Extend $a_n|_{V_n}$ to $X_0$ by affine extension on sub-triangles.  Call
the extension $g$.  Then $g$ takes the boundary of $X_0$ to a curve in $\Ce$
which contains $\frac{\sqrt{3}}{2}X_0$.
Since the boundary of $X_0$ is contractible
in $X_0$, it follows that $g(X_0)$ contains $\frac{\sqrt{3}}{2}X_0$,
and hence that
the area of $g(X_0)$ is strictly greater than $\frac{3}{4}$
times the area of $X_0$.

However, for any sub-triangle $S$ with barycenter $y$, if
$L_n(a_n) \le 1$ then
$g(S)$ is contained in the disk of radius $1/2^{n+1}$ about $a_n(y)$.
It is easy to see that
the area of $g(S)$ is thus at most $\frac{3}{4}$
times the area
of $S$.  Summing over $S$, we get that the area of $g(X_0)$ is at most
$\frac{3}{4}$
times the area of $X_0$.  This contradicts the previous paragraph, so we
conclude that $L_n(a_n) > 1$.
\end{proof}

Suppose that  $L_0(a)=L(a)=1$. Then there is an $f\in \Gamma$ with
$\limsup_{n\to \infty}L_n(f_{\frac{1}{n}})$ $=L_0(a)=1$.
Rescaling $f_{\frac{1}{n}}$ slightly, we may assume that
$L_n(f_{\frac{1}{n}})\le 1$ for all $n$.
Then Proposition~\ref{modified 4 point:prop} tells us that $\pa
f_{\frac{1}{n}}|_{V_n}-a|_{V_n}\pa \ge \frac{\sqrt{3}}{6} - \frac{1}{4}$.
This contradicts the assumption that $f\in \Gamma$.
Therefore $L$ does not coincide with
$L_0$ on the whole $C_{\Ce}(X_0)$.
\end{example}
\end{remark}


\section{Continuous fields of $C^*$-algebraic compact quantum \\metric
  spaces}
\label{COntFieldC*QCM:sec}

In this section we discuss the continuity of parameterized families of
$C^*$-algebraic compact quantum metric spaces with respect to $\dist_{\cq}$.
General criteria are given in Subsection~\ref{CriteriaC*GHD:sub}.
In Subsection~\ref{ContFieldC*QCMbyG:sub} we consider especially
families of $C^*$-algebraic compact quantum metric space
induced by
ergodic compact group actions.

\subsection{Criteria of $C^*$-algebraic quantum distance convergence}
\label{CriteriaC*GHD:sub}

The criteria we shall give are analogues of \cite[Theorems 1.2 and
7.1]{Li10}. A notion of continuous fields of compact quantum metric
spaces was introduced in \cite[Definition 6.4]{Li10} for families
of compact quantum metric spaces. In view of
Remark~\ref{Weaver:remark}, we give the following definition of
continuous fields of $C^*$-algebraic compact quantum metric
spaces. We refer the reader to \cite[Section 10.1 and 10.2]{Dixmier77} for basic
definitions and facts about continuous fields of $C^*$-algebras.

\begin{definition}  \label{cont field C*QCM:def}
Let $T$ be a locally compact Hausdorff space, and let $(\cA_t, L_t)$ be a
$C^*$-algebraic compact quantum metric space for
each $t\in T$. Let $\Gamma$ be the set of continuous sections of
 a continuous
field of $C^*$-algebras over $T$ with fibres $\cA_t$, and let
$\Gamma_{\sa}:=\{f\in \Gamma:f^*=f\}$ be the set of continuous sections
of  the associated continuous field of real Banach
spaces with fibres $(\cA_t)_{\sa}$.
 We call $(\{(\cA_t, L_t)\}, \Gamma)$ a
\emph{continuous field of $C^*$-algebraic compact quantum metric
  spaces} over $T$ if $(\{(A_t, L_t)\}, \Gamma_{\sa})$
is a continuous field of compact quantum metric spaces, \ie,
the unit section $t\mapsto e_{\cA_t}$ is in $\Gamma$, and for any $t'\in T, a\in A_{t'}$, and $\varepsilon>0$
there exists $f\in (\Gamma_{\sa})^L_{t'}$
such that $\pa f_{t'}-a\pa<\varepsilon$ and $L_{t'}(f_{t'})<L_{t'}(a)+\varepsilon$,
where
$$(\Gamma_{\sa})^L_{t'}=\{f\in \Gamma_{\sa}:
\mbox{ the function } t\mapsto L_t(f_t) \mbox{ is upper
semi-continuous at } t'\}.$$
\end{definition}

As an analogue of \cite[Theorem 1.2]{Li10}, we have the
following criterion for continuous fields of $C^*$-algebraic
compact quantum metric spaces to be continuous with respect to the
$C^*$-algebraic quantum distance:

\begin{theorem}  \label{criterion of C*conv:thm}
Let $(\{(\cA_t, L_t)\}, \Gamma)$ be a continuous field of
$C^*$-algebraic compact quantum metric spaces
over a locally compact Hausdorff space $T$.
Let $t_0\in T$, and let $\{f_n\}_{n\in \Ne}$ be a sequence in $\Gamma_{\sa}$ such that
$(f_n)_{t_0}\in \cD(A_{t_0})$ for each $n\in \Ne$ and
the set $\{(f_n)_{t_0}:n\in \Ne\}$ is dense in $\cD(A_{t_0})$.
Then the following are equivalent:

(1) $\dist_{\cq}(\cA_t, \cA_{t_0})\to 0$ as
$t\to t_0$;

(2)  $\dist_{\oq}(A_t, A_{t_0})\to 0$ as $t\to t_0$;

(3)  $\dist_{\GH}(\cD(A_t),
\cD(A_{t_0}))\to 0$ as
$t\to t_0$;

(4) for any $\varepsilon>0$, there is an $N$
such that the open $\varepsilon$-balls in $A_t$ centered at $(f_1)_t,
{\cdots}, (f_N)_t$ cover $\cD(A_t)$ for
all $t$ in some neighborhood $\mathcal{U}$ of $t_0$.
\end{theorem}
\begin{proof}
(1)$\Longrightarrow$(2) follows from Theorem~\ref{basic C*QGH:thm}.
(2)$\Longrightarrow$(3)$\Longrightarrow$(4)
follows from \cite[Theorem 1.2]{Li10}.
Hence we are left to show (4)$\Longrightarrow$(1).

As in the proof of \cite[Theorem 1.2]{Li10}
we may assume that $(f_n)_t\in \cD(A_t)$ for all $t\in T$ and $n\in \Ne$.
We can also find a normed space $V$ containing all $\cA_t$'s such
that for every $f\in \Gamma$ the map $t\mapsto f_t$ from $T$ to $V$
is continuous at
$t_0$. For any $n\in \Ne$ and $t\in T$
set $Y_{n,t}=\{(f_1)_t, {\cdots}, (f_n)_t\}$. Also
set $X_t=\cD(A_t)$ for all $t\in T$.

It is easy to see from (4) that there are a $R>0$ and a neighborhood
$\mathcal{U}_1$ of $t_0$ such that $R>r_{\cA_t}$ throughout
$\mathcal{U}_1$. Let $\varepsilon>0$ be given.
Pick $N$ and a neighborhood $\mathcal{U}\subseteq \mathcal{U}_1$
of $t_0$ for $\varepsilon$ as in (4).
Then $\dist^V_{\rH}(X_t, Y_{N, t})\le \varepsilon$
throughout $\mathcal{U}$. Let $t\in \mathcal{U}$. For any
$a, b\in X_t$ and $j, k\in \Ne$ we have
\begin{eqnarray*}
\pa ab-(f_j)_t(f_k)_t\pa &\le &\pa
a-(f_j)_t\pa \cdot \pa b\pa +\pa (f_j)_t\pa \cdot \pa b-(f_k)_t\pa \\
&\le &R( \pa
a-(f_j)_t\pa+\pa b-(f_k)_t\pa).
\end{eqnarray*}
Consequently
$$\dist^{V^{(3)}}_{\rH}((X_t)^m, (Y_{N, t})^m)\le \max(2R, 1)\dist^V_{\rH}(X_t,
  Y_{N, t})\le (2R+1)\varepsilon.$$
Clearly
$\dist^{V^{(3)}}_{\rH}((Y_{N,t})^m), (Y_{N, t_0})^m)\to 0$ as $t\to t_0$.
By shrinking $\mathcal{U}$ we may assume that $\dist^{V^{(3)}}_{\rH}((Y_{N,
      t})^m), (Y_{N, t_0})^m)<\varepsilon$ throughout $\mathcal{U}$.
Then
\begin{eqnarray*}
\dist_{\cq}(\cA_t, \cA_{t_0})
\le \dist^{V^{(3)}}_{\rH}((X_t)^m, (X_{t_0})^m)
\le (4R+3)\varepsilon
\end{eqnarray*}
for all $t\in \mathcal{U}$.
\end{proof}

It was pointed out in \cite[Section 11]{Li10} that the $\theta$-deformations
$M_{\theta}$
of Connes and Landi \cite{CL01} for a compact spin manifold $M$ form a natural continuous field
of $C^*$-algebras over the space of $n\times n$ skew-symmetric matrices $\theta$.
By \cite[Theorem 1.4]{Li10}, when $M$ is connected
the $M_{\theta}$'s form a continuous field of compact quantum metric spaces
and are continuous with respect to $\dist_{\oq}$.
Then Theorem~\ref{criterion of C*conv:thm} tells us that they are also continuous with respect to
$\dist_{\cq}$.

Combining \cite[Theorems 7.1]{Li10} and Theorems~\ref{basic C*QGH:thm} and \ref{criterion of C*conv:thm} together, we get

\begin{theorem}  \label{criterion 2 of C*conv:thm}
Let $(\{(\cA_t, L_t)\}, \Gamma)$ be a continuous field of
$C^*$-algebraic compact quantum metric spaces over a
locally compact Hausdorff space $T$. Suppose $R\ge r_{\cA_t}$ for
all $t\in T$. Let $t_0\in T$, and let $\{f_n\}_{n\in \Ne}$ be a sequence in
$\Gamma_{\sa}$ such that $(f_n)_{t_0}\in \cD_R(A_{t_0})$
for each $n\in \Ne$ and the set $\{(f_n)_{t_0}:n\in
\Ne\}$ is dense in $\cD_R(A_{t_0})$. Then the
following are equivalent:

(1) $\dist^R_{\cq}(\cA_t, \cA_{t_0})\to 0$ as $t\to
t_0$;

(2)  $\dist^R_{\oq}(A_t, A_{t_0})\to 0$ as $t\to t_0$;

(3) $\dist_{\GH}(\cD_R(A_t),
\cD_R(A_{t_0}))\to 0$ as
$t\to t_0$;

(4) for any $\varepsilon>0$, there is an $N$
such that the open $\varepsilon$-balls in $A_t$ centered at $(f_1)_t,
{\cdots}, (f_N)_t$ cover $\cD_R(A_t)$ for
all $t$ in some neighborhood $\mathcal{U}$ of $t_0$;

(5)  $\dist_{\cq}(\cA_t, \cA_{t_0})\to 0$ as $t\to
t_0$;

(6)  $\dist_{\oq}(A_t, A_{t_0})\to 0$ as $t\to t_0$.
\end{theorem}

In the theory of operator algebras usually continuous fields of $C^*$-algebras
are used to describe continuous families of algebras qualitatively.
Certainly, convergence under quantum distances is more concrete. So
one may expect the following:

\begin{proposition} \label{conv to cont field:prop}
Let $\{(\cA_n, L_n)\}_{n\in \Ne}$ be a sequence of $C^*$-algebraic
compact quantum metric spaces converging to
some $C^*$-algebraic compact quantum metric space
$(\cA_0, L_0)$ under $\dist_{\cq}$.
Then there is a
continuous field of $C^*$-algebraic compact quantum metric spaces
over $T=\{\frac{1}{n}:n\in \Ne\}\cup \{0\}$ with fibres
$(\cA_n, L_n)$ at $\frac{1}{n}$ and
$(\cA_0, L_0)$ at $0$.
\end{proposition}
\begin{proof}
We may assume that all the $L_n$'s and $L_0$ are closed.
Let $R>\sup\{r_{\cA_n}:n\in \Ne\}$. Then
$R\ge r_{\cA_0}$, and $\dist^R_{\cq}(\cA_n,
\cA_0)\to 0$ by Theorem~\ref{basic C*QGH:thm}.
For each $n$ let $V_n$ be a complex normed space containing both $\cA_n$ and
$\cA_0$ with
$\dist^{(V_n)^{(3)}}_{\rH}(\cD_R(A_n), \cD_R(A_0))$ $<\dist^R_{\cq}(\cA_n, \cA_0)+\frac{1}{n}$.
By \cite[Lemma 4.5]{Li10} we can find a complex normed space $V$
containing all these $V_n$'s with the copies of $\cA_0$
identified.
Then $$\dist^{V^{(3)}}_{\rH}(\cD_R(A_n), \cD_R(A_0))\to 0$$ as $n\to \infty$.
Let $\Gamma$ be the set of all continuous maps
$f$ from $T$ to $V$ such that $f(\frac{1}{n})\in \cA_n$ for
all $n$ and $f(0)\in \cA_0$. Clearly $\Gamma$ is the space of
continuous sections of a continuous field of Banach spaces over $T$
with fibres
$(\cA_n, L_n)$ at $\frac{1}{n}$ and
$(\cA_0, L_0)$ at $0$. As in Lemma~\ref{prod:lemma} it is easy
to see that for any $f, g\in \Gamma$ the sections $f\cdot g$ and $f^*$
are also in $\Gamma$. Thus this is a continuous field of
$C^*$-algebras. Also as in Lemma~\ref{lip:lemma} the restriction of
$L_0$ on $(\cA_0)_{\sa}$ is given by (\ref{limit L:eq}).
Thus this is a continuous field of $C^*$-algebraic compact quantum metric spaces.
\end{proof}

\begin{remark} \label{nuclear:remark}
\cite[Proposition 7.4]{Li10} is crucial in the proof of
(4)$\Longrightarrow$(1) in Theorem~\ref{criterion of C*conv:thm}.
As pointed out at the beginning of Section~\ref{C*QGH:sec}, one can
try to define a quantum distance, $\dist_{\rnu}$, as
$\inf\{\dist^{\cC}_{\rH}(h_{\cA}(\cD(A)), h_{\cB}(\cD(B)))\}$,
where the infimum is taken over all faithful $*$-homomorphisms $h_{\cA}$ and
$h_{\cB}$ of $\cA$ and $\cB$
into some
$C^*$-algebra $\cC$. Similarly, one defines
$\dist^R_{\rnu}$. One can check that $\dist_{\rnu}$ (resp. $\dist^R_{\rnu}$)
is also a metric on $\CCQM$ (resp. $\CCQM^R$)
 and that the forgetful map $(\CCQM,
\dist_{\rnu})\to (\CQM, \dist_{\q})$ is continuous (see
Theorem~\ref{basic C*QGH:thm}). Also $\dist_{\rnu}$ and
$\dist^R_{\rnu}$ define the same topology on $\CCQM^R$. The proofs
of  Proposition~\ref{dist_cq<dist_GH:prop} and Theorem~\ref{C*GH
GH:thm} hold with $\dist_{\cq}$ replaced by $\dist_{\rnu}$. But
without the $C^*$-algebraic analogue of \cite[Proposition
7.4]{Li10} we do not know whether (4) still implies (1) in
Theorem~\ref{criterion of C*conv:thm} with $\dist_{\cq}$ replaced
by $\dist_{\rnu}$. However, Blanchard has proved
\cite{Blanchard98} that every separable unital continuous field of
nuclear $C^*$-algebras over a compact metric space $T$ can be
subtrivialized.
Thus $\dist_{\rnu}$ behaves well at least for nuclear
$C^*$-algebras, and Theorems~\ref{criterion of C*conv:thm}
and \ref{criterion 2 of C*conv:thm}
hold with $\dist_{\cq}$ and $\dist^R_{\cq}$ replaced by $\dist_{\rnu}$
and $\dist^R_{\rnu}$
when $T$ is a compact metric space and each fibre $\cA_t$ is
nuclear. Limits of exact $C^*$-algebras under $\dist_{\rnu}$ are also exact:
\end{remark}

\begin{proposition} \label{exact:prop}
Let $\{(\cA_n, L_n)\}_{n\in \Ne}$ be a sequence of $C^*$-algebraic
compact quantum metric spaces converging to
some $C^*$-algebraic compact quantum metric space
$(\cA, L)$ under $\dist_{\rnu}$. If each $\cA_n$ is
exact, then so is $\cA$.
\end{proposition}
\begin{proof}
Let $R>\sup\{r_{\cA_n}:n\in \Ne\}$. Then
$R\ge r_{\cA}$, and $\dist^R_{\rnu}(\cA_n,
\cA_0)\to 0$.
For each $n$ let $\cC_n$ be a $C^*$-algebra containing both $\cA_n$ and
$\cA$ as $C^*$-subalgebras  with
$$\dist^{\cC_n}_{\rH}(\cD_R(A_n), \cD_R(A))<\dist^R_{\rnu}(\cA_n, \cA)+\frac{1}{n}.$$
Using amalgamation of $C^*$-algebras \cite[Theorem 3.1]{Blackadar80}
we can find a $C^*$-algebra $\cC$
containing all these $\cC_n$'s with the copies of $\cA$
identified.
Then $\dist^{\cC}_{\rH}(\cD_R(A_n), \cD_R(A))\to 0$ as $n\to \infty$. Taking a faithful $*$-representation
of $\cC$, we may assume that $\cC=B(H)$ for
some Hilbert space $H$. Denote by $1_H$ the identity operator on $H$.
Set $\tilde{\cA_n}=\Ce\cdot 1_H+\cA_n$, and let
$\tilde{\cA}=\Ce\cdot 1_H+\cA$.
Notice that $\tilde{\cA_n}=\Ce (1_H-e_{\cA_n})\oplus
\cA_n$. So $\tilde{\cA_n}$ is also exact.

Recall that a unital completely positive (u.c.p.) map
$\varphi:\cB_1\to \cB_2$ between unital
$C^*$-algebras is \emph{nuclear} if for any finite subset
$\mathcal{J}\subseteq \cB_1$ and any $\varepsilon>0$ there
are an integer $k$ and u.c.p. maps $\phi:\cB_1\to
M_k(\Ce), \, \psi:M_k(\Ce)\to \cB_2$ such that $\pa
(\psi\circ \phi)(x)-\varphi(x)\pa<\varepsilon$ for all $x\in
\mathcal{J}$. (For nonunital $C^*$-algebras, $\phi$ and $\psi$ are
required to be completely positive contractions.) A $C^*$-algebra
$\cB_1$ is \emph{nuclearly
  embeddable} if for some $C^*$-algebra $\cB_2$ there is a
faithful $*$-homomorphism $\varphi: \cB_1\to \cB_2$
such that $\varphi$ is nuclear. It is a theorem of Kirchberg \cite{Wassermann94} that a $C^*$-algebra
$\cB_1$ is exact if and only if $\cB_1$ is nuclearly
embeddable, and if and only if  any
faithful $*$-homomorphism $\varphi$ of $\cB_1$
into $B(H_1)$ of any Hilbert space $H_1$ is nuclear.
We will show that the inclusion $h_{\tilde{\cA}}:\tilde{\cA}\to B(H)$ is
nuclear. Then $\tilde{\cA}$ is exact. Since $\tilde{\cA}=\Ce (1_H-e_{\cA})\oplus
\cA$, the algebra $\cA$ is also exact.

Let a finite subset $\mathcal{J}\subseteq \tilde{\cA}$ and $\varepsilon>0$
be given. Since $A+iA+\Ce 1_H$ is dense in $\tilde{\cA}$ and
u.c.p. maps are contractions, we may assume that
$\mathcal{J}\subseteq A+\Re 1_H$. Notice that $\cD_R(A)$ is
absorbing in $A$, \ie $\Re_+\cD_R(A)=A$. Thus we may assume that $\mathcal{J}\subseteq
\cD_R(A)\cup \{1_H\}$. Take $N$ such that $\dist^{\cC}_{\rH}(\cD_R(A_N), \cD_R(A))<\frac{1}{3}\varepsilon$.
Then we can find a finite subset
$\mathcal{J}_N\subseteq \cD_R(A_N)\cup \{1_H\}$ with
$\dist^{\cC}_{\rH}(\mathcal{J}_N,\mathcal{J})<\frac{1}{3}\varepsilon$.
Since $\tilde{\cA_N}$ is exact, the inclusion $\tilde{\cA_N}\to B(H)$ is
nuclear. So there  are an
integer $k$ and u.c.p. maps
$\phi:\tilde{\cA_N}\to M_k(\Ce), \, \psi:M_k(\Ce)\to B(H)$
such that $\pa (\psi\circ \phi)(x)-\varphi(x)\pa<\frac{1}{3}\varepsilon$ for all
$x\in \mathcal{J}_N$. By Arveson's extension theorem \cite[Theorem 5.1.7]{ER00} we can extend
$\phi$ to a u.c.p. map from $B(H)$ to $M_k(\Ce)$ which we still denote by
$\phi$. Then it is easy to see that $\pa (\psi\circ \phi)(y)-\varphi(y)\pa<\varepsilon$ for all
$y\in \mathcal{J}$. Therefore the inclusion $h_{\tilde{\cA}}:\tilde{\cA}\to B(H)$ is
nuclear.
\end{proof}

\begin{example}[Quotient Field of a $C^*$-algebraic
Compact Quantum Metric Space]
\label{quotient field of C*QCM:eg}
Let $(\cA, L)$ be a closed $C^*$-algebraic compact quantum
metric space. Denote by $\mathcal{I}$
the set of all closed two-sided ideals not equal to $\cA$.
Let $(\{(B_t, L_t)\}, \Gamma)$ be the quotient field over $T$ for $(A, L|_A)$
  as constructed in
\cite[Example 6.3]{Li10}, where $T$ is the set of all
nonempty convex closed subsets of $S(A)=S(\cA)$ equipped with
the Hausdorff distance $\dist^{S(A)}_{\rH}$ and $(B_t, L_t)$ is the quotient
compact quantum metric space of $(A, L|_A)$ for each $t\in T$.
Denote by $\pi_t$ the quotient map $A\rightarrow B_t$ for each $t\in T$.
For any $I\in \mathcal{I}$
according to Lemma~\ref{quotient of C*QCM:lemma} the quotient
$\cA/I$ equipped with the quotient Lip-norm is a
closed $C^*$-algebraic compact quantum metric space.
 Now $S(\cA/I)$ is a convex
closed subset of $S(\cA)$. Clearly $I$ is determined by
$S(\cA/I)$ as
\begin{eqnarray*}
I=\{a\in \cA:a|_{S(\cA/I)}=0\}.
\end{eqnarray*}
So the map $I\mapsto S(\cA/I)$ from $\mathcal{I}$
to $T$
is injective.
In fact, the image of $\mathcal{I}$ is closed in
$(T, \dist^{S(\mathcal{A)}}_{\rH})$: let
$I_j$ be a sequence of in $\mathcal{I}$ such that $t_j=S(\cA/I_j)$
converges to some $t_0\in T$. Then $\pa \pi_{t_j}(a)\pa \to \pa
  \pi_{t_0}(a)\pa$ for all $a \in \cA_{\sa}$.
It follows immediately that the function $\pa \cdot \pa_{t_0}$ defined by
$\pa a\pa_{t_0}=(\pa \pi_{t_0}(aa^*)\pa)^{\frac{1}{2}}$ is a nontrivial
$C^*$-algebraic seminorm on $\cA$, \ie
$\pa a a'\pa_{t_0}\le \pa a\pa_{t_0} \pa a'\pa_{t_0}$
and $\pa aa^*\pa_{t_0}=\pa a\pa^2_{t_0}$
for all $a,a'\in \cA$. Denote by $I_0$ the kernel of $\pa \cdot
\pa_{t_0}$. Then $I_0$ is a closed two-sided ideal of
$\cA$. We claim that $t_0=S(\cA/I_0)$.
For any $\mu \in t_0$ clearly the evaluation at $\mu$
induces a state of $\cA/I_0$. Under the identification of
$S(\cA/I_0)$ with a closed convex subset of $S(\cA)$,
this state is just $\mu$ itself. Hence we see that $t_0\subseteq
S(\cA/I_0)$. For each $a\in \cA_{\sa}$ notice that
$\pa \pi_{t_0}(a)\pa=\pa a\pa_{t_0}=\pa
\pi_{S(\cA/I_0)}(a)\pa$.
By Lemma~\ref{cont to conv:lemma} below
we have
$t_0=S(\cA/I_0)$. Therefore
the image of $\mathcal{I}$ in $T$ is closed.
Identify $\mathcal{I}$ with its image, and let $(\{(B_t, L_t)\},
\Gamma_{\mathcal{I}})$ be the restriction of $(\{(B_t, L_t)\}, \Gamma)$ on
$\mathcal{I}$. Then $(B_{S(\cA/I)}, L_{S(\cA/I)})$ is
  the associated quantum metric space of $(\cA/I,
  L_{\cA/I})$.
Notice that for any $a\in
\cA$ the map $I\mapsto \pa \pi_{\cA/I}(a)\pa =
(\pa \pi_{\cA/I}(aa^*)\pa)^{\frac{1}{2}}$ is continuous.
Then it is easy to see that $(\{(\cA/I, L_{\cA/I})\},
\Gamma_{\mathcal{I}}+i\Gamma_{\mathcal{I}})$ is a continuous field of
$C^*$-algebraic compact quantum metric spaces.
 We call it the {\it quotient field of $(\cA, L)$}.
It was noticed in
\cite[Example 6.3]{Li10} that
$\dist_{\q}(B_{S(\cA/I)}, B_{S(\cA/I_0)})\to 0$ as
$I\to I_0$. By \cite[Theorems 1.1]{Li10}
and Theorem~\ref{criterion of C*conv:thm} we see that $\dist_{\cq}( \cA/I,
\cA/I_0)\to 0$ as $I\to I_0$.
\end{example}

\begin{remark} \label{Fell:remark}
The construction of the continuous field of order-unit spaces in
\cite[Example 6.3]{Li10} and that of the continuous field of
$C^*$-algebras in example~\ref{quotient field of C*QCM:eg} work in
general. For any compact Hausdorff space $X$, the set $\SUB(X)$ of closed
nonempty subsets of $X$ is a compact Hausdorff space equipped with
the Hausdorff topology, which is defined using the pseudometrics
$\dist^{\rho}_{\rH}$ for all the continuous pseudometrics $\rho$ on $X$ (see
\cite{Fell62} for a generalization).
For any order-unit space $(A, e)$ (see for example \cite[Section 2]{Li10}),
it is clear that
the set of all nonempty closed convex subsets of the state-space $S(A)$, which we
still denote by $T$ as in \cite[Example 6.3]{Li10},
is a closed subspace of $\SUB(S(A))$.
So $T$ is a compact
Hausdorff space with the relative topology.
As in \cite[Example 6.3]{Li10}
the restrictions of
$a\in A$ to closed convex subsets of $S(A)$ still generate a
continuous field of order-unit spaces over $T$, which we denote
by $(\{B_t\}, \Gamma)$. For a unital $C^*$-algebra $\cA$
the argument in
Example~\ref{quotient field of C*QCM:eg} still works to show that the set
of all proper closed two-sided ideals $\mathcal{I}$
is a closed subset of the set $T$ of all nonempty convex closed subsets of
$S(\cA_{\sa})=S(\cA)$, and the quotient
images for $a\in \cA$ generate a continuous field of
$C^*$-algebras over $\mathcal{I}$. When $\cA$ is nonunital,
let $\tilde{\cA}=\cA+\Ce e$ be $\cA$ with
a unit adjoined. Then $\mathcal{I}(\tilde{\cA})=
\mathcal{I}(\cA)\cup \{A\}$. So $\mathcal{I}(\cA)$
becomes a locally compact Hausdorff space with $\infty$ being the
ideal $\cA$ itself, and we still get a continuous field of
$C^*$-algebras over $\mathcal{I}(\cA)$ generated by the
quotient images of $a\in \cA$. In both cases, the natural
$*$-homomorphism from $\cA$ to the algebras of
global sections vanishing at $\infty$ is injective and the composed
homomorphism to each fibre algebra is surjective. This means that
$\cA$ is represented as a \emph{full algebra of operator
  fields} in the sense of \cite[page 236]{Fell61}. In fact, one can check easily
that our topology
on $\mathcal{I}$ coincides with that on the space of
norm-functions (with $0$ removed) in \cite[page 243]{Fell61}.
Now it is clear that the usage of the
quantum metric on $A$ (resp. $\cA$)
in \cite[Example 6.3]{Li10}
(resp. Example~\ref{quotient field of C*QCM:eg}) is just to make $T$
(resp. $\mathcal{I}$) into a metric space and to endow each fibre space
of the continuous field of order-unit spaces (resp. $C^*$-algebras)
with a quantum metric (resp. $C^*$-algebraic quantum metric).
\end{remark}

\begin{lemma} \label{cont to conv:lemma}
Let $(A, e)$ be an order-unit space,
and let $\{\mathfrak{X}_n\}$
be a net of closed convex subsets of $S(A)$. Let $\pi_n:A\rightarrow
\Af_{\Re}(\mathfrak{X}_n)$ be the restriction map, where $\Af_{\Re}(\mathfrak{X}_n)$
is the space of continuous affine $\Re$-valued functions on $\mathfrak{X}_n$
and we identify $A$ with a dense subspace of $\Af_{\Re}(S(A))$ via the pairing between
$A$ and $S(A)$.
If $\pa a\pa=\lim_{n\to
  \infty}\pa \pi_n(a)\pa$ for all $a\in A$, then
$\mathfrak{X}_n$ converges to $S(A)$ under
the Hausdorff topology in $\SUB(S(A))$
when $n\to \infty$.
\end{lemma}
\begin{proof} Suppose that $\mathfrak{X}_n$ does not converge to $S(A)$
under
the Hausdorff topology
when $n\to \infty$. Then there are some
$\varepsilon>0$ and a subnet of $\{\mathfrak{X}_n\}$
converging to some closed
convex $\mathfrak{X}\subsetneq S(A)$. Without loss of generality we may assume that
this subnet is $\{\mathfrak{X}_n\}$ itself.
Pick $\mu\in S(A)\setminus
\mathfrak{X}$. By the Hahn-Banach theorem there exist $a'\in
\bar{A}$ and $t'\in \Re$ with $\mu(a')>t'\ge \nu(a')$ for all
$\nu\in \mathfrak{X}$. Since $A$ is dense in $\bar{A}$, we may also
find $a\in A$ and $t\in \Re$ such that $\mu(a)>t\ge \nu(a)$ for all
$\nu\in \mathfrak{X}$. Then $\frac{\mu(a)+t}{2}\ge \nu(a)$ for all
$\nu\in \mathfrak{X}_n$ when $n$ is big enough. Consequently,
 $\mu(a-s\cdot e)>\frac{\mu(a)+t}{2}-s\ge \nu(a-s\cdot e)\ge0$ for all
$\nu\in \mathfrak{X}_n$, where
$s=\min \{a(\nu):\nu\in S(A)\}$.
Therefore $\pa a-s\cdot e\pa\ge \mu(a-s\cdot
e)>\frac{\mu(a)+t}{2}-s\ge \pa \pi_n(a-s\cdot e)\pa$, which
contradicts the assumption that $\lim_{n\to \infty}\pa a-s\cdot e\pa=\pa \pi_n(a-s\cdot
e)\pa$. So $\mathfrak{X}_n$ does converge to $S(A)$ under
the Hausdorff topology
when $n\to \infty$.
\end{proof}

\begin{example}[Berezin-Toeplitz Quantization] \label{BT:eg}
Let $(M, \omega)$ be a K\"ahler manifold \cite{Cannas01, GH94},
\ie $M$ is a
complex manifold and $\omega$ is a K\"ahler form. This means that
$\omega$ is a positive, non-degenerate closed $2$-form of type
$(1,1)$. Locally, if $\dim_{\Ce}(M)=d$ and $z_1, {\cdots}, z_d$ are
local holomorphic coordinates,  then $\omega$  can be written as
\begin{eqnarray*}
\omega=i\sum^d_{j, k=1}g_{jk}(z)dz_j\wedge d\bar{z}_k,&  g_{jk}\in
C^{\infty}(M, \Ce),
\end{eqnarray*}
where the matrix $(g_{jk}(z))$ is a positive
definite hermitian matrix for each $z$.  Consider
triples $(E, h, \nabla)$, with $E$ a holomorphic line
bundle, $h$ a hermitian metric on $E$ , and $\nabla$ a
connection which is compatible with the metric and the complex
structure. With respect to local holomorphic coordinates and a local
holomorphic frame $s$ of the bundle, this means that
$\nabla=\partial+\partial \log\hat{h}+\bar{\partial}$, where
$\hat{h}(z)=h_z(s(z), s(z))$. The curvature of $\nabla$ is defined as
\begin{eqnarray*}
F(X,Y)=\nabla_X\nabla_Y-
\nabla_Y\nabla_X-\nabla_{[X, Y]}.
\end{eqnarray*}
The K\"ahler manifold $(M, \omega)$ is called
\emph{quantizable} \cite[Section 2]{Schlichenmaier98}
if there is such a triple $(E, h, \nabla)$ with
\begin{eqnarray} \label{quantizable:eq}
F=-i\omega.
\end{eqnarray}
This includes all the compact Riemann surfaces and the complex
projective spaces. Now let $M$ be a quantizable compact  K\"ahler
manifold, and let $\Omega=\frac{1}{d!}\omega^d$ be the volume form on $M$.
Then the space of holomorphic sections $\Gamma_{\hol}(M, E)$ is a finite
dimensional complex vector space \cite[page 152]{GH94}.
On the space of smooth sections $\Gamma_{\infty}(M, E)$ we have the
inner product
\begin{eqnarray*}
 \left< \varphi, \psi\right>:=\int_Mh(\varphi, \psi)\Omega.
\end{eqnarray*}
Denote by $L^2(M, E)$ the $L^2$-completion of $\Gamma_{\infty}(M, E)$, and
denote by
$\Pi$ the orthogonal projection $L^2(M, E)\rightarrow \Gamma_{\hol}(M, E)$. For any
$f\in C_{\Ce}(M)$ the \emph{Toeplitz operator} $T_f$ is defined as
\begin{eqnarray*}
T_f:=\Pi(f\cdot ): \Gamma_{\hol}(M, E)\rightarrow \Gamma_{\hol}(M, E).
\end{eqnarray*}
It is easy to see that $f\mapsto T_f$ is a unital completely positive map from $C_{\Ce}(M)$ to
$B(\Gamma_{\hol}(M, E))$, the space of bounded operators on
$\Gamma_{\hol}(M, E)$. In particular, $T_{f^*}=(T_f)^*$ for all
$f\in C_{\Ce}(M)$. We may replace $E$ by the
$n$-th tensor powers $E^n:=E^{\otimes n}$ and apply above construction for every $n$.
In this way we get linear maps
\begin{eqnarray*}
T^{(n)}:C_{\Ce}(M)\rightarrow B(\Gamma_{\hol}(M, E^n)), & f\mapsto T^{(n)}_f.
\end{eqnarray*}
The condition (\ref{quantizable:eq}) implies that $E$ is ample. Replacing
$E$ and $\omega$ by $E^{n_0}$ and $n_0\omega$ respectively
for some $n_0$ we may assume that $E$ is very ample.
Then $T^{(n)}$ is surjective for all $n\in \Ne$ \cite[Proposition 4.2]{Schlichenmaier94}.
For any metric $\rho_M$ on $M$ inducing the
manifold topology, we have the associated $C^*$-algebraic
compact quantum metric space $(C_{\Ce}(M), L)$.
Then for  any $n\in \Ne$ the $C^*$-algebra $B(\Gamma_{\hol}(M,
E^n))$, equipped with the quotient Lip-norm $L_{1/n}$ constructed in
Lemma~\ref{quotient of C*QCM:lemma}, is a $C^*$-algebraic quantum
compact metric space.
One has $\pa T^{(n)}(f)\pa \to \pa f\pa$ \cite[Theorem 4.1]{Schlichenmaier94} and $\pa
(T^{(n)}(f))\cdot (T^{(n)}(g))-T^{(n)}(f\cdot g)\pa\to 0$ \cite[page 291]{Schlichenmaier94}
\cite[Theorem 3]{Schlichenmaier98}
as $n\to \infty$
for any $f, g\in C_{\Ce}(M)$. Consequently,
the sections given by the images of $f$ under
$T^{(n)}$ for $f\in C_{\Ce}(M)$ generate a continuous field of
$C^*$-algebras over $T'=\{\frac{1}{n}:n\in \Ne\}\cup
\{0\}$, with fibres $B(\Gamma_{\hol}(M, E^n))$ at $1/n$ and $C_{\Ce}(M)$
at $0$. In fact, it was shown in \cite{Schlichenmaier94}
 that this is a
strict quantization of the symplectic structure on $M$, though we do not
need this fact here.
For each $f\in C_{\Re}(M)$, we have
$L_{\frac{1}{n}}(T^{(n)}(f))\le L(f)$. So this is a
continuous field of $C^*$-algebraic compact quantum metric spaces. If we identify
$S(B(\Gamma_{\hol}(M, E^n)))$ with the corresponding convex closed subset in $S(C_{\Ce}(M))$,
then Lemma~\ref{cont to conv:lemma} tells us that $S(B(\Gamma_{\hol}(M, E^n)))$ converges to $S(C_{\Ce}(M))$ under
$\dist^{S(C_{\Ce}(M))}_{\rH}$, and this is a special case of \cite[Example 6.10]{Li10}.
By \cite[Theorem 1.1]{Li10} and Theorem~\ref{criterion of C*conv:thm}
 we get that $\dist_{\cq}(B(\Gamma_{\hol}(M, E^n)), C_{\Ce}(M))\to 0$ as $n\to \infty$.
Let us mention that when $M$ is an integral coadjoint orbit
of a compact connected semisimple Lie group with the standard
symplectic form, this becomes the Berezin quantization in \cite{Rieffel01}
(see also \cite[Example 10.12]{Li10}).
For such case, the quotient Lip-norms $L_{1/n}$
and the induced metrics on the state-spaces $S(B(H_n))$
are also considered in \cite{ZS00}. Note that the quotient Lip-norm here
is different from the Lip-norm induced by the ergodic action as in
\cite{Rieffel01}.
\end{example}

\subsection{Continuous fields of $C^*$-algebraic compact quantum
  metric spaces induced by ergodic compact group actions}
\label{ContFieldC*QCMbyG:sub}

Let $G$ be a compact group with a fixed length function $\mathnormal{l}$,
\ie, a continuous
real-valued function, $\mathnormal{l}$, on $G$ such that
\begin{eqnarray*}
\mathnormal{l}(xy)&\le &\mathnormal{l}(x)+\mathnormal{l}(y)\mbox{ for
  all } x, y\in G \\
\mathnormal{l}(x^{-1})&=& \mathnormal{l}(x) \mbox{ for all }x\in G \\
\mathnormal{l}(x)&=& 0 \mbox{ if and only if } x=e_G,
\end{eqnarray*}
where $e_G$ is the identity of $G$.

Let $\cA$ be a unital
$C^*$-algebra, and let $\alpha$ be a strongly continuous ergodic action of $G$ on
$\cA$ by automorphisms.
In \cite[Theorem 2.3]{Rieffel98b} Rieffel showed
that $\cA$ is a $C^*$-algebraic compact quantum metric space equipped with
the seminorm $L$ defined by
\begin{eqnarray} \label{def of L:eq}
L(a)=\sup \{\frac{\pa \alpha_x(a)-a\pa}{\mathnormal{l}(x)}: x\in G, x\neq e_G\}.
\end{eqnarray}

Let $(\{\cA_t\}, \Gamma)$ be a continuous field
of $C^*$-algebras over a locally compact Hausdorff space $T$,
and let $\alpha_t$ be a strongly continuous action of $G$ on $\cA_t$
for each $t\in T$.
Recall that $\{\alpha_t\}$ is a \emph{continuous field
of strongly continuous actions} of $G$ on $(\{\cA_t\}, \Gamma)$
if the actions $\{\alpha_t\}$ give rise to a strongly continuous
action of $G$ on $\Gamma_{\infty}$ \cite[Definition 3.1]{Rieffel89b},
where $\Gamma_{\infty}$ is the space
of continuous sections vanishing at $\infty$.
If each $\cA_t$ is unital and
$\alpha_t$ is ergodic, we say that this is a \emph{field of ergodic
actions}. In such case, one can show easily that the unit section is continuous.

Denote by $\hat{G}$ the dual group of $G$.
For any strongly continuous
action of $G$ on a $C^*$-algebra $\cA$ and $\gamma \in \hat{G}$,
denote by
$\mul(\cA, \gamma)$ the multiplicity of $\gamma$ in $\cA$.
Combining \cite[Theorem 1.3]{Li10} and
Theorem~\ref{criterion of C*conv:thm} together, we have

\begin{theorem} \label{criterion of cont field of C*action:thm}
Let $\{\alpha_t\}$ be a continuous field of strongly continuous ergodic
actions of $G$ on a continuous field of unital C$^*$-algebras
$(\{\cA_t\}, \Gamma)$
over a locally compact Hausdorff space $T$.
Then the induced field $(\{(\cA_t, L_t)\}, \Gamma)$ is a continuous field
of $C^*$-algebraic compact quantum metric spaces.
For any  $t_0\in T$  the following are equivalent:

(1) $\lim_{t\to t_0}\mul(\cA_t, \gamma)=
\mul(\cA_{t_0}, \gamma)$ for all $\gamma \in \hat{G}$;

(2) $\limsup_{t\to t_0}\mul(\cA_t, \gamma)\le
\mul(\cA_{t_0}, \gamma)$ for all $\gamma \in \hat{G}$;

(3) $\dist_{\oq}(A_t, A_{t_0})\to 0$ as $t\to t_0$;

(4) $\dist_{\cq}(\cA_t, \cA_0)\to 0$ as $t\to t_0$.
\end{theorem}

In particular, Rieffel's continuity of noncommutative tori \cite[Theorem 9.2]{Rieffel00}
and matrices converging to integral coadjoint orbits of compact connected semisimple Lie groups
\cite[Theorem 3.2]{Rieffel01} hold with respect to $\dist_{\cq}$ (see \cite[Examples 10.11 and 10.12]{Li10}).

As a consequence of Theorem~\ref{criterion of cont field of
  C*action:thm}
 we show that when $T$ is a compact metric space,
the map $t\mapsto (\cA_t, L_t)$ from $T$ to $\CCQM$
is continuous at most points of $T$:

\begin{corollary} \label{2nd cat:coro}
Let $\{\alpha_t\}$ be a continuous field of strongly continuous ergodic
actions of $G$ on a continuous field of unital C$^*$-algebras
$(\{\cA_t\}, \Gamma)$
over a  compact metric space $T$. Then there is
a nowhere dense $F_{\sigma}$ subset $Z\subseteq T$
such that for any $t_0\in T\setminus Z$,
$\dist_{\cq}(\cA_t, \cA_{t_0})\to 0$ as $t\to t_0$,
\end{corollary}
\begin{proof}
By \cite[Proposition 2.1]{HLS81}, for each $\gamma \in \hat{G}$ the
multiplicity of $\alpha_t$ in $\cA_t$ is no bigger than
$\dim(\gamma)$. Set $X_{\gamma, \dim(\gamma)}=\{t \in T:
\mul(\cA_t, \gamma)=\dim(\gamma)\}$. Since the multiplicity
function is lower semi-continuous over $T$ by \cite[Lemma 10.6]{Li10},
$X_{\gamma, \dim(\gamma)}$ is an open subset of $T$.
Inductively, for $1\le j\le \dim(\gamma)$ set $X_{\gamma,
\dim(\gamma)-j}=\{t \in T\setminus
\overline{\cup^{j-1}_{k=0}X_{\gamma, \dim(\gamma)-k}} :
\mul(\cA_t, \gamma)=\dim(\gamma)-j\}$. Then the $X_{\gamma,
\dim(\gamma)-j}$'s are all open subsets of $T$. Set
$X_{\gamma}=\cup_{0\le j \le \dim{\gamma}}X_{\gamma,
\dim(\gamma)-j}$. Then $X_{\gamma}$ is the subset of $T$ with
constant multiplicity with respect to $\gamma$ locally. Clearly
$X_{\gamma}$ is open, and its complement $Z_{\gamma}$
has no interior points.
By \cite[Remark 8.1]{Li10}, $G$ is metrizable and
hence $L^2(G)$ is separable.
Since every $\gamma\in \hat{G}$ appears in the left regular
 representation, $\hat{G}$ is countable.
 Denote by $Z$ the union of
all these $Z_{\gamma}$'s.
Then $Z$ satisfies the
requirement according to Theorem~\ref{criterion of cont field of
  C*action:thm}.
\end{proof}


\begin{thebibliography}{10}

\bibitem{Alfsen71}
E. M. Alfsen. {\it Compact Convex Sets and Boundary Integrals.}
Ergebnisse der Mathematik und ihrer Grenzgebiete, Band 57. Springer-Verlag, New York-Heidelberg, 1971.

\bibitem{Blackadar80}
B. E. Blackadar. Weak expectations and nuclear $C^*$-algebras.
{\it Indiana Univ. Math. J.}  {\bf 27} (1978), no. 6, 1021--1026.

\bibitem{BC91}
B. E. Blackadar and J. Cuntz.  Differential Banach algebra norms and smooth subalgebras of $C^*$-algebras.
{\it J. Operator Theory}  {\bf 26} (1991), no. 2, 255--282.

\bibitem{Blanchard98}
E. Blanchard. Subtriviality of continuous fields of nuclear $C^*$-algebras.
{\it J. Reine Angew. Math.}  {\bf 489} (1997), 133--149. math.OA/0012128.

\bibitem{Schlichenmaier94}
M. Bordemann, E. Meinrenken, and M. Schlichenmaier. Toeplitz quantization of K\"ahler manifolds and ${\rm gl}(N)$, $N\to \infty$
  limits.
{\it Comm. Math. Phys.}  {\bf 165} (1994), no. 2, 281--296.
hep-th/9309134.

\bibitem{BBI01}
D. Burago, Y. Burago, and S. Ivanov. {\it A Course in Metric Geometry}.
Graduate Studies in Mathematics, 33.
\newblock American Mathemetical Society, Providence, RI, 2001.

\bibitem{Cannas01}
A. Cannas~da Silva. {\it Lectures on Symplectic
Geometry.}
Lecture Notes in Mathematics, 1764.
Springer-Verlag, Berlin, 2001.

\bibitem{CC83}
M.-D. Choi and E. Christensen. Completely order isomorphic and close $C^*$-algebras need not be $*$-isomorphic.
{\it Bull. London Math. Soc.}  {\bf 15} (1983), no. 6, 604--610.

\bibitem{Connes89}
A. Connes.  Compact metric spaces, Fredholm modules, and hyperfiniteness.
{\it Ergodic Theory Dynam. Systems}  {\bf 9} (1989), no. 2, 207--220.

\bibitem{Connes94}
A. Connes. {\it Noncommutative Geometry}.
Academic Press, Inc., San Diego, CA, 1994.

\bibitem{CL01}
A. Connes and G. Landi. Noncommutative manifolds, the instanton algebra and isospectral deformations.
{\it Comm. Math. Phys.}  {\bf 221} (2001), no. 1, 141--159.
math.QA/0011194.

\bibitem{Conway90}
J. B. Conway. {\it A Course in Functional Analysis.}
Second edition. Graduate Texts in Mathematics, 96. Springer-Verlag, New York, 1990.

\bibitem{Dixmier77}
J. Dixmier. {\it $C^*$-algebras}.
Translated from the French by Francis Jellett. North-Holland Mathematical Library, Vol. 15.
North-Holland Publishing Co., Amsterdam-New York-Oxford, 1977.

\bibitem{ER00}
E. G. Effros and Z.-J. Ruan. {\it Operator Spaces.}
London Mathematical Society Monographs. New Series, 23.
The Clarendon Press, Oxford University Press, New York, 2000.

\bibitem{Fell61}
J. M. G. Fell. The structure of algebras of operator fields.
{\it Acta Math.}  {\bf 106} (1961), 233--280.

\bibitem{Fell62}
J. M. G. Fell. A Hausdorff topology for the closed
  subsets of a locally compact non-Hausdorff space.
{\it Proc. Amer. Math. Soc.}  {\bf 13} (1962), 472--476.

\bibitem{GH94}
P. Griffiths and J. Harris. {\it Principles of Algebraic Geometry}.
Reprint of the 1978 original. Wiley Classics Library. John Wiley \& Sons, Inc., New York, 1994.

\bibitem{Gromov81}
M. Gromov.  Groups of polynomial growth and expanding maps.
{\it Inst. Hautes \'Etudes Sci. Publ. Math.}  {\bf 53} (1981), 53--73.

\bibitem{Gromov99}
M. Gromov. {\it Metric Structures for Riemannian and non-Riemannian Spaces.}
Based on the 1981 French original. With appendices by M. Katz, P. Pansu and S. Semmes.
Translated from the French by Sean Michael Bates. Progress in Mathematics, 152. Birkh\"auser Boston, Inc., Boston, MA, 1999.


\bibitem{HLS81}
R. H\o egh-Krohn, M. B. Landstad, and  E. St\o rmer.  Compact ergodic groups of automorphisms.
{\it Ann. of Math. (2)}  {\bf 114} (1981), no. 1, 75--86.

\bibitem{Johnson94}
B. E. Johnson.  Near inclusions for subhomogeneous $C^*$-algebras.
{\it Proc. London Math. Soc. (3)}  {\bf 68} (1994), no. 2, 399--422.

\bibitem{Kadison51}
R. V. Kadison. A representation theory for commutative topological algebra.
{\it Mem. Amer. Math. Soc.} (1951), no. 7.


\bibitem{KS72}
R. V. Kadison and D. Kastler. Perturbations of von Neumann algebras. I. Stability of type.
{\it Amer. J. Math.}  {\bf 94} (1972), 38--54.


\bibitem{Kerr02}
D. Kerr. Matricial quantum Gromov-Hausdorff distance.
{\it J. Funct. Anal.}  {\bf 205}  (2003),  no. 1, 132--167. math.OA/0207282.

\bibitem{KL}
D. Kerr and H. Li. On Gromov-Hausdorff convergence for operator
metric spaces.  math.OA/0411157.



\bibitem{Li10}
H. Li. Order-unit quantum Gromov-Hausdorff distance.
math.OA/0312001 v2.

\bibitem{Li9}
H. Li. $\theta$-deformations as compact quantum metric spaces.
math.OA/0311500 v2.

\bibitem{Rieffel03O}
N. Ozawa and M. A. Rieffel. Hyperbolic group $C^*$-algebras and free-product $C^*$-algebras as compact quantum metric spaces.
{\it Canad. J. Math,} to appear. math.OA/0302310.

\bibitem{PR81}
J. Phillips and I. Raeburn.  Perturbations of
  $C^*$-algebras. II.
{\it Proc. London Math. Soc. (3)}  {\bf 43} (1981), no. 1, 46--72.

\bibitem{Rieffel89b}
M. A. Rieffel. Continuous fields of $C^*$-algebras coming from group cocycles and actions.
{\it Math. Ann.}  {\bf 283}  (1989),  no. 4, 631--643.

\bibitem{Rieffel98b}
M. A. Rieffel. Metrics on states from actions of
compact groups.
{\it Doc. Math.} {\bf 3}  (1998), 215--229 (electronic).
math.OA/9807084.


\bibitem{Rieffel99b}
M. A. Rieffel. Metrics on state spaces.
{\it Doc. Math.}  {\bf 4} (1999), 559--600 (electronic).
math.OA/9906151.


\bibitem{Rieffel00}
M. A. Rieffel. Gromov-Hausdorff distance for
quantum  metric spaces.
{\it Mem. Amer. Math. Soc.}  {\bf 168}  (2004),  no. 796, 1--65.
 math.OA/0011063.

\bibitem{Rieffel01}
M. A. Rieffel. Matrix algebras converge to the
sphere for  quantum Gromov-Hausdorff distance.
{\it Mem. Amer. Math. Soc.}  {\bf 168}  (2004),  no. 796, 67--91.
math.OA/0108005.

\bibitem{Rieffel02}
M. A. Rieffel. Group $C^*$-algebras as compact quantum metric spaces.
{\it Doc. Math.}  {\bf 7}  (2002), 605--651 (electronic). math.OA/0205195.


\bibitem{Rieffel03}
M. A. Rieffel. Compact quantum metric spaces.
math.OA/0308207.


\bibitem{Sakai96}
T. Sakai. {\it Riemannian Geometry}.
Translated from the 1992 Japanese original by the author. Translations of Mathematical Monographs, 149.
American Mathematical Society, Providence, RI, 1996.


\bibitem{Schlichenmaier98}
M. Schlichenmaier. Berezin-Toeplitz quantization of compact K\"ahler manifolds.
In: {\it Quantization, Coherent States, and Poisson Structures(Bia\l
  owie\.za, 1995)},   101--115, PWN, Warsaw, 1998.
q-alg/9601016.

\bibitem{Wassermann94}
S. Wassermann. {\it Exact $C^*$-algebras and Related Topics}.
Lecture Notes Series, 19. Seoul National University, Research Institute of Mathematics, Global Analysis Research Center, Seoul, 1994.


\bibitem{Weaver99}
N. Weaver. {\it Lipschitz Algebras}.
World Scientific Publishing Co., Inc., River Edge, NJ, 1999.

\bibitem{Weaver01}
N. Weaver. Personal communications.
(2001).

\bibitem{ZS00}
K. \.Zyczkowski and W.  S\l omczy\'nski.  The Monge metric on the sphere and geometry of quantum states.
{\it J. Phys. A}  {\bf 34} (2001), no. 34, 6689--6722.

\end{thebibliography}
\end{document}